\documentclass[twoside,11pt]{article}

\usepackage{blindtext}

\usepackage{jmlr2e}

\usepackage{smile}
\usepackage{amsmath}
\usepackage{dsfont}
\usepackage{mathtools}
\usepackage{bbm}
\usepackage{bm}
\usepackage{enumerate}
\usepackage{caption}
\usepackage{longtable}
\usepackage[linesnumbered,ruled,vlined]{algorithm2e}
\usepackage{color}

\mathtoolsset{showonlyrefs}

\newif\ifReviseMode\ReviseModefalse
\newcommand{\revise}[1]{\ifReviseMode{\color{blue} #1}\else{ #1}\fi}

\newcommand{\eqd}{\stackrel{\text{d}}{=}}
\newcommand{\supp}{\textrm{supp}}

\newcommand{\dataset}{{\cal D}}

\usepackage{lastpage}
\jmlrheading{24}{2023}{1-\pageref{LastPage}}{9/22; Revised
6/23}{9/23}{22-1032}{Ethan X. Fang, Yajun Mei, Yuyang Shi, Qunzhi Xu and Tuo Zhao}

\ShortHeadings{Pivotal Estimation of Linear Discriminant Analysis in High Dimensions}{Fang, Mei, Shi, Xu and Zhao}
\firstpageno{1}

\begin{document}

\title{Pivotal Estimation of Linear Discriminant Analysis\\ in High Dimensions}

\author{\name Ethan X. Fang$^\dagger$
\email ethan.fang@duke.edu 
       \AND
      \name Yajun Mei$^\ddagger$
      \email ymei3@gatech.edu 
       \AND
      \name Yuyang Shi$^\ddagger$
      \email yyshi@gatech.edu 
      \AND
      \name Qunzhi Xu$^\ddagger$
      \email xuqunzhi@gatech.edu
      \AND
      \name Tuo Zhao$^\ddagger$
      \email tourzhao@gatech.edu\\
      \addr 
      $\dagger$ \addr Department of Biostatistics and Bioinformatics, Duke University\\
      $\ddagger$ \addr School of Industrial and Systems Engineering,
      Georgia Tech
      }
\editor{Mladen Kolar}

\maketitle

\begin{abstract}
We consider the linear discriminant analysis problem in the high-dimensional settings. In this work, we propose PANDA(\underline{P}ivot\underline{A}l li\underline{N}ear \underline{D}iscriminant \underline{A}nalysis), a tuning-insensitive method  in the sense that it requires very little
effort to tune the  parameters. Moreover, we prove that PANDA achieves the optimal convergence rate in terms of
both the estimation error and misclassification rate.
Our theoretical results are backed up by thorough numerical studies using both simulated and real datasets. In comparison with the existing methods, we observe that our proposed PANDA yields equal or better performance, and requires substantially less effort in parameter tuning. 
\end{abstract}



\section{Introduction}
We consider the linear discriminant analysis problem with  $n_0$ samples $(X_i^{(0)})_{i=1}^{n_0}$ from class 0 and $n_1$ samples $(X_i^{(1)})_{i=1}^{n_1}$ from class 1. In particular, consider the Gaussian case where $X_i^{(\ell)}\sim N(\mu^{(\ell)},\Sigma), \ell=0,1.$    Under the ideal setting where all parameters $\mu^{(0)}, \mu^{(1)}, \Sigma$ are pre-specified, the Bayes rule classifies a new sample $Z$ by
\begin{align*}
f^*(Z) = \mathbbm{1}\big\{(Z-\mu_m)^\top\Sigma^{-1}\mu_d>0\big\},
\end{align*}
where $\mu_m = (\mu^{(0)} + \mu^{(1)})/2$ and $\mu_d = (\mu^{(1)} - \mu^{(0)}),$ and is proved to be optimal in terms of misclassification rate, see \cite{anderson1962introduction}. However, the Bayes rule is often not practical, as in reality the parameters are always unknown and  need to be estimated.

Under the classical low-dimensional setting $p < n$, we estimate $\mu^{(0)}$, $\mu^{(1)}$ and $\Sigma^{-1}$ by their sample versions, and use the plug-in  Bayes rule to classify the new sample. In particular, let $\widehat{\mu}^{(\ell)}$'s and $\widehat{\Sigma}$ be  the the sample means and the pooled sample covariance matrix, and let $\widehat{\mu}_m = (\widehat{\mu}^{(0)} + \widehat{\mu}^{(1)})/2$, $\widehat{\mu}_d = (\widehat{\mu}^{(1)} - \widehat{\mu}^{(0)})$. Given a new sample $Z$, the following rule
\begin{equation*}
\widehat{f}(Z) = \mathds{1}\left\{\hat{\mu}_d^\top\hat{\Sigma}^{-1}\left(Z-\hat{\mu}_m\right)>0\right\},
\end{equation*}
asymptotically achieves the optimal Bayesian risk. Unfortunately, this method is inapplicable to high-dimensional settings where $p \gg n$ because it is difficult to estimate $\Sigma^{-1}$ due to the singularity of $\widehat{\Sigma}$. Such high dimensionality issues  exist unavoidably in many critical modern scenarios such as genomics, and  it is important to develop efficient methods for LDA in high dimensions.

Several methods have been developed in the literature for high-dimensional LDA with sparsity assumptions imposed, which are common in many real-world applications such as the fMRI decoding and biomarker identification \citep{yamashita2008sparse,shi2009sparse}.
The existing methods can be further divided into two tracks based on the different sparsity assumptions.
The first track is to assume that   $\Sigma$ is sparse and estimate $\mu_d = \mu^{(1)}-\mu^{(0)}$ and $\Sigma$ separately. A simple approach is the naive Bayes rule or independence rule discussed in \cite{bickel2004some}. \cite{tibshirani2002diagnosis}, and \cite{fan2008high} proposed the nearest shrunken centroid method and the Features Annealed Independence Rules (FAIR) respectively for selecting significant features. Also see the sparse linear discriminant analysis (SLDA) proposed in \cite{shao2011sparse}.


Another track of work assumes the sparsity of the discriminant direction $\beta^* = \Sigma^{-1}\mu_d$ and directly estimates $\beta^*$ from the samples. 
\cite{witten2011penalized} and \cite{clemmensen2011sparse} proposed the sparse discriminant analysis method with multiple classes by imposing fused LASSO penalty and elastic net penalty respectively.
\cite{mai2012direct} proposed to estimate $\beta^*$ by minimizing an $\ell_1$-penalized least square loss, and
\cite{fan2012road} proposed the regularized optimal affine discriminant (ROAD) method. 

Existing theoretical results in the literature of high-dimensional LDA often require
the knowledge of unknown population. For the better understanding, here we present the linear programming discriminant (LPD) rule in \cite{cai2011direct} with more details.
The LPD rule provides an estimator  $\widehat{\beta}$ for $\beta^*$  by solving the following linear optimization problem 
\begin{align*}
\widehat{\beta} \in \mathop{\arg\min}_{\beta\in\mathbb{R}^p}\ \norm{\beta}_1,\quad \text{subject~to}\quad\norm{\widehat{\Sigma}\beta - \widehat{\mu}_d}_\infty\leq \lambda\widehat{\sigma}_{\max},
\end{align*}
with  $\hat{\sigma}_{\max} = \sqrt{\max_j\hat{\Sigma}_{jj}}$ and tuning parameter  $\lambda$. The authors show that to ensure the fast 
convergence rate of $\widehat\beta,$ a reasonable choice of $\lambda$ would be
\begin{equation*}
\lambda = O\left(\Delta\sqrt{\frac{\log p}{n}}\right),
\end{equation*} 
where $\Delta=\sqrt{\beta^{*\top}\Sigma\beta^*}.$
In practice, this choice of $\lambda$ heavily relies on the unknown population quantity~$\Delta$, which takes substantial effort to tune. 
 To reduce the tuning effort, \cite{tony2019high} propose the adaptive linear discriminant analysis (AdaLDA) rule, which is a two-stage method that achieves the minimax optimal convergence rate in both the estimation error and misclassification rate. Specifically, the AdaLDA rule solves a two-stage problem: in the first stage it constructs an estimator $\hat\Delta$ for $\Delta$ and in the second stage the estimator is plugging into the LPD framework to obtain the estimator for $\beta^*$.
 


In this paper, we propose a novel one-stage method for high-dimensional linear discriminant analysis named \textbf{PANDA} (\underline{P}ivot\underline{A}l li\underline{N}ear \underline{D}iscriminant \underline{A}nalysis). Our method is tuning-insensitive, in the sense that it automatically adapts to the population pattern and requires less effort to tune. Motivated by \cite{gautier2011high}
for high-dimensional linear regression, 
the proposed PANDA method simultaneously estimates $\beta^*$ and $\Delta$  by solving a single convex optimization problem, and is shown to
attain the same minimax optimal convergence rate  as  the AdaLDA method. Moreover, our detailed numerical results show that the PANDA method achieves similar or more competitive performance than the LPD and AdaLDA methods in terms of $\beta^*$ estimation error and misclassification rate, with less cost of computational time.

It is worth mentioning that the topic of variable selection  has also been investigated in high-dimensional LDA. For example, \cite{kolar2014optimal} established the optimal results of variable selection for sparse discriminant analysis in \cite{mai2012direct} and the ROAD estimator in \cite{fan2012road}, and 
\cite{Gaynanova2014Optimal} further extended the result to the multi-group sparse discriminant analysis. We also include some numerical studies investigating the variable selection properties of PANDA in Section \ref{numerical_results}. 

\vspace{5pt}

\noindent{\bf Paper Organization.} The rest of this paper is organized as follows. In Section \ref{background}, we briefly review the LDA problem  and the AdaLDA rule. In Section \ref{method}, we propose the PANDA method. In Section~\ref{theory}, we provide theoretical justifications of  PANDA.
In Section \ref{numerical_results}, we present the numerical studies. 
In Section \ref{extension K-class}, we discuss the extension of our PANDA method to the multiple-class LDA problem.
In Section \ref{proof}, we provide proofs of our main results. 
We conclude the paper in Section \ref{discussion}.

\vspace{5pt}

\noindent\textbf{Notations.} Let $v= (v_1, \cdots, v_p)^\top\in \mathbb{R}^p$ be a $p-$dimensional real vector. We define the following vector norms: $\norm{v}_1 = \sum_{j=1}^p |v_j|$, $\norm{v}^2_2 = \sum_{j=1}^p v_j^2$, and  $\norm{v}_\infty = \max_{1\leq j\leq p}|v_j|$.  For $p\in \mathbb{N}$, we denote by $[p]$ the set $\left\{1,2,\cdots,p \right\}$. For $j\in [p]$, let $e_j$ be the $j-$th canonical basis in $\mathbb{R}^p$. For $S\subseteq[p]$, let $v_S$ denote the the subvector of $v$ confined to $S$, and $S^c$ denotes the complement of $S$. For a matrix $\Sigma\in\RR^{p\times p}$, $\Sigma \succ 0$  denotes that $\Sigma$ is symmetric and positive definite, and  $\lambda_{\min}(\Sigma)$ and $\lambda_{\max}(\Sigma)$  denote the smallest and the largest eigenvalue of $\Sigma$, respectively. We let $\bm{0}$ and $\bm{1}$  denote  vectors with all the entries equal to $0$ and $1$, respectively. We use $\mathds{1}\{\cdot\}$ to denote the indicator function.

\section{Background}\label{background}
In this section, we provide necessary mathematical background.
For better presentation, we split this section into two subsections. We review the problem setup of LDA in Section \ref{problem setup}, and the AdaLDA method in Section \ref{AdaLDA}. 

\subsection{Problem Setup}\label{problem setup}

We consider the linear discriminant analysis problem with  $n_0$ samples $(X_i^{(0)})_{i=1}^{n_0}$ from class 0 and $n_1$ samples $(X_i^{(1)})_{i=1}^{n_1}$ from class 1. In particular, consider the Gaussian case where $X_i^{(\ell)}\sim N(\mu^{(\ell)},\Sigma), \ell=0,1.$  
Our goal is to find a linear discriminant rule $f_{\alpha,\beta}(\cdot)$ such that given a new sample $Z$, we predict the class label of $Z$ by
\begin{equation*}
f_{\alpha,\beta}(Z) = \mathds{1}\left\{
\beta^\top (Z-\alpha) > 0
\right\},
\end{equation*}
with some $\alpha,\beta\in\RR^{p}$. For simplicity, we assume the two classes  have equal prior weights, i.e., $\PP(\textrm{$Z$ is from Class $0$})=\PP(\textrm{$Z$ is from Class $1$})=1/2.$
Then the misclassification rate of $f_{\alpha,\beta}(\cdot)$ can be written as
\begin{align}
\mathcal{R}(f_{\alpha,\beta}) &= \frac{1}{2} \mathbb{P}_{Z\sim N(\mu^{(0)}, \Sigma)}({f}_{\alpha,\beta}(Z) = 1) + \frac{1}{2} \mathbb{P}_{Z\sim N(\mu^{(1)}, \Sigma)}({f}_{\alpha,\beta}(Z) = 0)\notag\\
& = \frac{1}{2} \Phi \left(-\frac{{\beta}^\top\left(\alpha-\mu^{(0)}\right) }{\sqrt{{\beta}^\top\Sigma{\beta}}} \right) + \frac{1}{2} \Phi \left(-\frac{{\beta}^\top\left(\mu^{(1)} - \alpha\right) }{\sqrt{{\beta}^\top\Sigma{\beta}}} \right), \label{risk_formula}
\end{align}
where $\Phi$ is the CDF of the standard Gaussian distribution.

The optimal misclassification rate (also known as the Bayes error) is achieved by the Fisher's discriminant rule $f_{\alpha^*,\beta^*}(\cdot)$ with $\alpha^* = (\mu^{(0)}+\mu^{(1)})/{2}$ and $\beta^*=\Sigma^{-1}\left(\mu^{(1)}-\mu^{(0)}\right).$
Accordingly, the optimal misclassification rate is $\mathcal{R}^* = \Phi(-\Delta/2)$, where $\Delta = \sqrt{\beta^{*\top}\Sigma\beta^*} = \sqrt{\mu_d^\top\Sigma^{-1}\mu_d}$ is the signal-noise ratio of the classification problem.

\subsection{The AdaLDA method}\label{AdaLDA}


In this subsection, we review the AdaLDA method proposed in \cite{tony2019high}, which is tuning-insensitive and serves as a good comparison to our method.
Let the sample means  and the pooled covariance matrix be
\begin{align*}
\hat{\mu}^{(\ell)}=\frac{1}{n_\ell}\sum_{i=1}^{n_\ell}X_{i}^{(\ell)}\quad\textrm{and}\quad\hat{\Sigma}=\frac{1}{n_0+n_1}\sum_{\ell = 0,1}\sum_{i=1}^{n_\ell}(X_{i}^{(\ell)}-\hat{\mu}^{(\ell)})(X_{i}^{(\ell)}-\hat{\mu}^{(\ell)})^\top
.
\end{align*}
The AdaLDA method estimates $\beta^*$ through two stages. In the first stage, AdaLDA solves the following linear optimization problem to obtain an initial estimator $\tilde{\beta}$,
\begin{align} \label{AdaLDA program 1}
\begin{split}
\tilde{\beta} \in \mathop{\arg\min}_{\beta}\quad&\norm{\beta}_1,\\
\text{subject~to}\quad &\norm{\hat{\Sigma}\beta - \hat{\mu}_d}_{\infty} \leq 4\hat{\sigma}_{\max}\cdot\sqrt{\frac{\log p}{n}}\cdot\left(\lambda\beta^\top\hat{\mu}_d\ + 1\right),
\end{split}
\end{align}
where $n=\min(n_0, n_1)$, $\lambda$ is a tuning parameter, $\hat{\mu}_d = \hat{\mu}^{(1)}-\hat{\mu}^{(0)}$ is the difference of the sample means, and $\hat{\sigma}_{\max} = \sqrt{\max_j\hat{\Sigma}_{jj}}$. The initial esstimator $\tilde{\beta}$ is used to construct an estimator $\hat{\Delta} = \sqrt{|\tilde{\beta}^{\top}\hat{\mu}_d|}$ for $\Delta$. In the second stage, AdaLDA solves another linear optimization problem to obtain the final estimator $\hat\beta$ 
\begin{align*}
\hat{\beta} \in \mathop{\arg\min}_{\beta}\quad&\norm{\beta}_1,\notag\\
\text{subject~to}\quad&|e_j^\top (\hat{\Sigma}\beta - \hat{\mu}_d)|\leq 4\hat{\sigma}_{\max}\cdot\sqrt{\frac{\log p}{n}}\cdot \sqrt{\lambda\hat{\Delta}^2 + 1},~ \text{ for all }j\in[p].
\end{align*}
With  $\hat{\beta}$ and $\hat{\mu}_m = \left(\hat{\mu}^{(0)}+\hat{\mu}^{(1)}\right)/2$, AdaLDA constructs the linear discriminant rule
$f_{\hat{\mu}_m, \hat{\beta}}$.

With a slight abuse of the notation, we let $\cR(\hat{\beta}) = \mathcal{R}(f_{\hat{\mu}_m,\hat{\beta}})$. 
Since the tuning parameters in the two steps do not depend on any unknown population quantities, the AdaLDA method is tuning-insensitive. 
Assuming $\beta^*$ contains at most~$s$ nonzero entries, \cite{tony2019high} prove that  under some mild assumptions, 
by choosing  $\lambda$
 as a proper constant, both $\hat{\beta}$ and $\mathcal{R}(\hat{\beta})$ achieve the minimax optimal rates of convergence that
\begin{align*}
\norm{\hat{\beta} - \beta^*}_2 =\cO_{P}\left(\Delta\sqrt{\frac{s\log p}{n}}\right)\quad\textrm{and}\quad
\mathcal{R}(\hat{\beta})- \mathcal{R}^*= \cO_P\left(\exp\left(-\frac{\Delta^2}{8}\right)\Delta \frac{s\log p}{n}\right).
\end{align*}

\section{The PANDA Method}\label{method}

In this section, we propose PANDA, a one-stage and tuning-insensitive method  for linear discriminant analysis in high dimensions. To begin with, we would like to first recall the LPD method \cite{cai2011direct}, which motivates our formulation.
 Specifically, the LPD method estimates $\beta^*$ by solving the following linear optimization problem that
\begin{align} \label{LPD}
\mathrm{\textbf{LPD}}:\quad \widehat{\beta}\in \mathop{\arg\min}_{\beta\in \mathbb{R}^p}~\norm{\beta}_1,\quad\text{subject~to~}\quad \norm{\widehat{\Sigma}\beta-\widehat{\mu}_d}_\infty\leq \lambda\hat{\sigma}_{\max}.  
\end{align} 
As discussed earlier, the tuning parameter $\lambda$ in \eqref{LPD} depends on  the unknown population quantity $\Delta = \sqrt{\beta^{*\top}\Sigma\beta^*}$, which is difficult to tune in practice.

To address this issue, we introduce $\tau$ as an estimator of $\Delta$, and plug it into \eqref{LPD}, as inspired by the pivotal method for high-dimensional linear regression in \cite{gautier2011high}. This leads to the following optimization problem 
\begin{align}\label{equality_constraint_program}
(\widehat{\beta},\widehat{\tau})\in \mathop{\arg\min}_{\beta\in \mathbb{R}^p, \tau\in\mathbb{R}}~\norm{\beta}_1,
\quad\text{subject~to~}\quad\norm{\widehat{\Sigma}\beta-\widehat{\mu}_d}_\infty\leq \lambda\hat\sigma_{\max}(\tau+1),\quad \sqrt{\beta^\top \widehat{\Sigma}\beta}=\tau.
\end{align}
The optimization problem in \eqref{equality_constraint_program} is nonconvex due to the  quadratic equality constraint $\sqrt{\beta^\top \widehat{\Sigma}\beta}=\tau$. Thus, we propose to relax the equality constraint into an inequality constraint, and obtain
\begin{align}
\label{inequality_constraint_program}
(\widehat{\beta},\widehat{\tau})\in \mathop{\arg\min}_{\beta\in \mathbb{R}^p, \tau\in\mathbb{R}}~\norm{\beta}_1,\quad\text{subject~to~}\quad \norm{\widehat{\Sigma}\beta-\widehat{\mu}_d}_\infty\leq \lambda \hat\sigma_{\max}(\tau+1),\quad \sqrt{\beta^\top \widehat{\Sigma}\beta}\leq\tau. 
\end{align}
However, as the objective function in \eqref{inequality_constraint_program}  is free of $\tau$, $\tau$ can be arbitrarily large. In fact, \eqref{inequality_constraint_program} admits a trivial solution $\widehat{\beta} = {0}$ when  $\widehat{\tau}$ is larger than $\lambda^{-1}\norm{\hat{\mu}_d}_{\infty}-1$, which makes \eqref{inequality_constraint_program}  inapplicable.

To solve this problem, we introduce an additional penalty term $c\tau^2$ to the objective  in \eqref{inequality_constraint_program}, which leads to the following PANDA's  formulation:
\begin{align} \label{con:problem}
\mathrm{\textbf{PANDA}}:\quad(\hat{\beta},\hat{\tau}) \in \mathop{\arg\min}_{\beta\in\mathbb{R}^p,\tau\in\mathbb{R}}\quad&\norm{\beta}_1+c\tau^2,\notag\\
\quad\textrm{subject~to}\quad&\norm{\widehat{\Sigma}\beta-\hat{\mu}_d}_{\infty}\leq\lambda\hat{\sigma}_{\max}(\tau+1), \quad\sqrt{\beta^\top \widehat{\Sigma}\beta}\leq \tau, 
\end{align}
where $c>0$ and $\lambda>0$ are two tuning parameters. 
Note that different from the linear penalty term used in Gautier's pivotal method, our penalty term is quadratic in $\tau.$ In fact, we can show that to guarantee the tuning-insensitivity of our PANDA method, the penalty term on $\tau$ must be quadratic. We provide more detailed discussion in
Section \ref{penalty_term} of the supplementary material.



Note that both our proposed PANDA method and the AdaLDA method adopt the similar idea of plugging in an estimator of the unknown quantity $\Delta$ to the tuning parameter $\lambda$ in the LPD method to achieve tuning-insensitivity. The main difference is that AdaLDA constructs the estimator for $\Delta$ in a separate linear program \eqref{AdaLDA program 1}, while PANDA estimates $\beta^*$ and $\Delta$ in a single convex program.




We point out that the problem in \eqref{con:problem} is a second order conic optimization problem. By introducing auxiliary variables $w\in\mathbb{R}^p$ and $u\in \mathbb{R}$, the problem in \eqref{con:problem} is equivalent to the following optimization problem 
\begin{align}
\min_{\beta, \tau, w, u}\quad& \sum_{j=1}^p w_j + cu,\label{socp}\\
\text{subject~to}\quad & -w_j\leq \beta_j\leq w_j,\quad -\lambda\widehat{\sigma}_{\max}(\tau + 1)\bm{1}\leq \widehat{\Sigma}\beta-\widehat{\mu}_d\leq \lambda\widehat{\sigma}_{\max} (\tau + 1)\bm{1},\notag\\
&\norm{\widehat{\Sigma}^{1/2}\beta}_2\leq \tau,\quad\sqrt{\tau^2 + \frac{1}{4}(1-u)^2}\leq \frac{1}{2}(1+u).\notag
\end{align}
Such a second order conic optimization problem is convex, and can be solved in a polynomial time using the interior point method \citep{nesterov1994interior}. Computationally, we also provide an efficient scheme in Algorithm \ref{ADMM} using the  alternating direction method of multipliers (ADMM) following \cite{boyd2011distributed} to solve \eqref{socp}.  We provide more details on the derivation of the algorithm  in Section~\ref{ADMM_scheme} of the supplementary material.
\begin{algorithm}
	\caption{ADMM with proximal method for solving problem \eqref{socp}}
	\begin{algorithmic}\label{ADMM}
		\STATE{\textbf{Input:} Sample mean difference $\widehat{\mu} = \widehat{\mu}^{(1)} - \widehat{\mu}^{(0)}$; Pooled sample covariance matrix $\widehat{\Sigma}$; Tuning parameter $c,~\lambda$; Initialization $\beta^0,~\tau^0~u^0,~v^0,~w^0,~s^0$;
		Penalty parameter $\rho>0$; Primal step size $\eta>0$; Number of iterations $T$.}\\
		\FOR{$t=1,2,\cdots,T$}
		\STATE $\beta^t \leftarrow \beta^{t-1} -\eta\nabla_{\beta} L_\rho(\beta^{t-1},u^{t-1},v^{t-1},w^{t-1},\tau^{t-1},s^{t-1})$
		\STATE $u^t\leftarrow \Pi_{\mathcal{C}_1}[
		u^{t-1} - \eta\nabla_u L_\rho(\beta^t,u^{t-1},v^{t-1},w^{t-1},\tau^{t-1},s^{t-1})]$
		\STATE $v^t\leftarrow \Pi_{\mathcal{C}_2}[v^{t-1} - \eta\nabla_v L_\rho(\beta^t,u^t,v^{t-1},w^{t-1},\tau^{t-1},s^{t-1})]$
		\STATE $\tilde{\tau}^{t} \leftarrow \tau^{t-1} -\eta\nabla_{\tau}L_{\rho}(\beta^t,u^{t},v^{t},w^{t-1},\tau^{t-1},s^{t-1})$
		\STATE $\tilde{w}^{t}\leftarrow w^{t-1} - \eta\nabla_w L_\rho(\beta^t,u^t,v^t,w^{t-1},\tau^{t-1},s^{t-1})$
		\STATE $(w^t,\tau^t)\leftarrow \Pi_{\mathcal{C}_2}(\tilde{w}^{t},\tilde{\tau}^{t})$
		\STATE $s^{t}\leftarrow s^{t-1} + A_\beta \beta^{t} + A_u u^t + A_v v^t + A_w w^t  + A_\tau \tau^t - b^t$
		\ENDFOR\\
  \STATE \textbf{Output } $(\beta^T,\tau^T,w^T,u^T)$
	\end{algorithmic}
\end{algorithm}

\section{Statistical Properties}\label{theory}
In this section, we establish theoretical guarantees for our proposed PANDA method. For notational simplicity, we denote
\begin{align*}
\mu_m &= (\mu^{(0)} + \mu^{(1)})/2,\quad\widehat{\mu}_m = (\widehat{\mu}^{(0)} + \widehat{\mu}^{(1)})/2,\quad\mu_d = \mu^{(1)} - \mu^{(0)},\\
\widehat{\mu}_d &= \widehat{\mu}^{(1)} - \widehat{\mu}^{(0)},\quad\sigma_{\max} = \max_j(\Sigma_{jj})^{1/2},\quad\widehat{\sigma}_{\max} = \max_j(\widehat{\Sigma}_{jj})^{1/2}.
\end{align*}
Without loss of generality, here we only consider the case where $n_0 = n_1 = n$, and our results can be easily extended to the general case where $n_0 \neq n_1$.  We  require the following weak sparsity condition on $\beta^*:$  
\begin{align}\label{l_q ball definition}
\beta^* \in \mathbb{B}_q(R)\coloneqq \left\{\beta\in \mathbb{R}^p: \sum_j |\beta_j|^q\leq R\right\},
\end{align}
where  $q\in[0,1)$ and $R$ can scale with $n$ and $p$. Note that when $q=0$, $\mathbb{B}_q(R)$ is reduced to the class of $R$-sparse vectors, i.e., $\mathbb{B}_0(R)\coloneqq \left\{\beta\in \mathbb{R}^p: \sum_j \mathds{1}\{\beta_j\neq 0\}\leq R\right\}$.
We also need to  impose the following two mild assumptions.
\begin{assumption}
	\label{identifiablity_assumption}
	There exists a  constant $a$ such that $\norm{\mu_d}_\infty\geq a>0$.
\end{assumption}

\begin{assumption}\label{Sigma_eigenvalue_assumption}
	There exists some $M$ such that $M^{-1} \leq \lambda_{\min}(\Sigma) \leq \lambda_{\max}(\Sigma)\leq M$.
\end{assumption}
Essentially, Assumption \ref{identifiablity_assumption} requires the two classes to be distinguishable, and
Assumption \ref{Sigma_eigenvalue_assumption} requires the covariance matrix $\Sigma$ to be sufficiently well-conditioned, as its condition number is upper bounded by $M^2$.

We are now ready to present the theoretical guarantees of the PANDA method in \eqref{con:problem}. Let us begin with
the convergence rates of $\hat{\beta}$ and $\hat{\tau}$.

\begin{theorem}[Parameter Estimation]\label{one_norm_two_norm_theorem}
	Suppose that Assumption
	 \ref{Sigma_eigenvalue_assumption} hold, and $\beta^*\in \mathbb{B}_q(R)$ for some $q\in[0,1)$ and some  $R>0$. Let $\left(\widehat{\beta},\widehat{\tau} \right)$ be an optimal solution of \eqref{con:problem}. Given
	\begin{align} \label{parameter_setting}
	c = \frac{1}{8\left(\norm{\widehat{\mu}_d}_\infty + 4\widehat{\sigma}_{\max}\sqrt{\frac{\log p}{n}}\right)}, 
	\quad \lambda = 20 \sqrt{\frac{\log p}{n}},
	\end{align} 
	for sufficiently large $n$ such that 
	\begin{align}\label{sample_complexity_1_norm}
	n\geq C^{(1)}\cdot a^{-2}\Delta^{2}\sigma_{\max}^2M^{2+\frac{1}{1-q}}R^{\frac{2}{1-q}}\log p
	\end{align}
	where $C^{(1)}$ is an absolute constant, we have, with probability goes to 1,
	\begin{subequations}
	\begin{align}
\norm{\hat{\beta}-\beta^*}_1 &\leq C_1 \cdot (\Delta+1)\left(\sigma_{\max} M\right)^{1-q}R\left(\frac{\log p}{n} \right)^{(1-q)/2},\label{1-norm_bound}\\
	\norm{\hat{\beta}-\beta^*}_2 &\leq C_2\cdot (\Delta+1) \left(\sigma_{\max} M\right)^{1-q/2}\sqrt{R}\left(\frac{\log p}{n}\right)^{1/2-q/4},\label{2-norm_bound}\\
	\frac{|\widehat{\tau}^2 -\Delta^2|}{\Delta^2} &\leq C_3\cdot (1+\Delta^{-1})\sigma_{\max}^{1-q/2}M^{3/2-q}R\Big(\frac{\log p}{n}\Big)^{(1-q)/2}, \label{Delta_bound}
	\end{align}
	\end{subequations}
	 where $C_1$, $C_2$ and $C_3$ are  positive constants.
	\label{1-norm_theorem}
\end{theorem}

	Note that our proposed PANDA method is tuning-insensitive, as the chosen tuning parameters $c$ and $\lambda$ in \eqref{parameter_setting} do not depend on any unknown population quantity. 
In the next theorem, we show that the sample complexity requirement \eqref{sample_complexity_1_norm} can be relaxed under some more restrictive conditions.

	\begin{theorem}\label{one_norm_two_norm_theorem_revise}
		Suppose that Assumption
		\ref{Sigma_eigenvalue_assumption} holds, and $\beta^*\in \mathbb{B}_q(R)$ for some $q\in[0,1)$ and some  $R>0$. Let $\left(\widehat{\beta},\widehat{\tau} \right)$ be an optimal solution to problem \eqref{con:problem}. When $\widehat{\tau} = \sqrt{\widehat{\beta}^\top\widehat{\Sigma}\widehat{\beta}}$, given
		\begin{align*}
		c = \frac{1}{8\left(\norm{\widehat{\mu}_d}_\infty + 4\widehat{\sigma}_{\max}\sqrt{\frac{\log p}{n}}\right)}, 
		\quad \lambda = 20 \sqrt{\frac{\log p}{n}},
		\end{align*} 
		for sufficiently large $n$ such that 
		\begin{align}\label{sample_complexity_revise}
		n\geq C^{(2)}\cdot a^{-2}\Delta^{2}\sigma_{\max}^2M^{2+\frac{1}{1-q}}R^{\frac{2}{2-q}}\log p,
		\end{align}
		where $C^{(2)}$ is an absolute constant, we have, with probability goes to 1,
		\begin{subequations}
			\begin{align}
			\norm{\widehat{\beta} - \beta^*}_2&\leq C_1\cdot (\Delta+1)(\sigma_{\max}M)^{1-q/2}\sqrt{R}\left(\frac{\log p}{n}\right)^{1/2-q/4},\label{2-norm_bound_revise}\\
			\frac{|\widehat{\tau}^2 -\Delta^2|}{\Delta^2} &\leq C_2\cdot (1+\Delta^{-1})\sigma_{\max}^{1-q/2}M^{3/2-q}\sqrt{R}\Big(\frac{\log p}{n}\Big)^{(1-q)/2}, \label{Delta_bound_revise}
			\end{align}
		\end{subequations}
		where $C_1$ and $C_2$ are  positive constants.
		\label{1-norm_theorem_revise}
	\end{theorem}

\textbf{}Note that in the above theorem, we impose the additional assumption that $\widehat{\tau} = \sqrt{\widehat{\beta}^\top \widehat{\Sigma}\widehat{\beta}}$, i.e. the second inequality constraint of PANDA is active at the optimal solution. We point out that in practice, we can numerically verify if this assumption indeed holds. Also, in our later simulations, we find that this assumption holds when the tuning parameters are properly chosen.



We next compare our results with \cite{tony2019high} for $q=0$. Note that \cite{tony2019high} consider the following parameter space of $\beta^*$ and $\Sigma$,
\begin{align}\label{AdaLDA parameter space}
\Theta_s = \big\{&(\beta^*,\Sigma):\beta^*\in \mathbb{R}^p,~\Sigma\in\mathbb{R}^{p\times p},~|\supp(\beta^*)|\leq s,\notag\\
&M^{-1}\leq\lambda_{\min}(\Sigma)\leq\lambda_{\max}(\Sigma)\leq M, \Delta\geq c_L >0 \big\},
\end{align}
where $M$ and $c_L$ are absolute constants that do not scale with $n$, $p$ and $s$. 
They then establish the following minimax lower bound,
\begin{align*}
\inf_{\widehat{\beta}}\sup_{(\beta^*,\Sigma)\in \Theta_s}\mathbb{E}\left[\norm{\widehat{\beta}-\beta^*}_2 \right]\geq C_M\cdot \Delta\sqrt{\frac{s\log p}{n}},
\end{align*}
where the infimum is taken over any estimator $\widehat{\beta}$ based on the samples, and $C_M$ is some constant depending on $M$. Under such a setting, both AdaLDA and PANDA are minimax optimal in terms of $\beta^*$ estimation.
When $M$ is allowed to scale with $n$, $p$ and $s$, the PANDA method still attains the same rates of convergence for parameter estimation as the AdaLDA method.
Specifically, we follow  the same analysis in \cite{tony2019high} and rewrite their results with explicit dependence on $M$ as follows,
	\begin{align*}
		\norm{\widehat{\beta}-\beta^*}_2& = \cO_P\left(\sigma_{\max} M\Delta\sqrt{\frac{s\log p}{n}}\right),\\
		\frac{|\widehat{\Delta}^2-\Delta^2|}{\Delta^2}& =
		\cO_P\left(\sigma_{\max}M^{3/2}
		\sqrt{\frac{s\log p}{n}}\right).
	\end{align*}
In addition, to ensure the above rates of convergence with high probability, the sample size $n$ needs to satisfy that
\begin{align*}
n = \cO_P\left(\sigma_{\max}^2 M^3 \Delta^2s\log p\right).
\end{align*}
As can be seen, in Theorem \ref{one_norm_two_norm_theorem_revise}, our convergence rates \eqref{2-norm_bound_revise} and \eqref{Delta_bound_revise} matches the convergence rates in \cite{tony2019high} with the same order of sample complexity.


Next, let us establish an upper bound for the misclassification rate of the obtained estimator $\hat{\beta}$ in the PANDA method.
\begin{theorem}[Misclassification Rate]	\label{risk_bound_theorem}
	 Under the identical conditions as in Theorem \ref{one_norm_two_norm_theorem} or \ref{one_norm_two_norm_theorem_revise}, we have, with probability goes to 1,
\begin{align*}
\mathcal{R}(\hat{\beta}) - \mathcal{R}^*\leq 
C\cdot \exp\left(-\frac{\Delta^2}{8}\right) \sigma_{\max}^{-q}M^{3-q}\Delta
 R \left(\frac{\log p}{n}\right)^{1-q/2}
\end{align*}
where $C$ is an absolute positive constants.
\end{theorem}
	When $q=0$ and $R=s$, \cite{tony2019high} consider the parameter space of $\beta^*$ and $\Sigma$ defined in \eqref{AdaLDA parameter space}, where $M$ is a constant,
	and establish the following minimax lower bound
	\begin{align*}
	\inf_{\hat{f}} \sup_{(\beta^*,\Sigma)\in \Theta_s}
	\cR(\hat{f}) - \cR^* \geq C\cdot \exp\left(-\frac{\Delta^2}{8}\right)\Delta^{-1}\frac{s\log p}{n},
	\end{align*}
	where the infimum is taken over any linear discriminant rule $\hat{f}$ based on the samples. Under such a setting, both AdaLDA and PANDA attain the minimax optimal rates of convergence for the misclassification rate that
	\begin{align*}
	\cR(\widehat{\beta})- \cR^* = \cO_P\left(\exp\left(-\frac{\Delta^2}{8}\right)
	M^3\Delta \frac{s\log p}{n} \right).
	\end{align*}

\begin{remark}
The probability of the convergence rates in Theorems \ref{one_norm_two_norm_theorem}, \ref{one_norm_two_norm_theorem_revise}   and \ref{risk_bound_theorem} is due to the uncertainty of data, which is addressed in Lemma \ref{events_lemma} and Lemma \ref{Gaussian design RE condition lemma} in later analysis. As a summary, the probability of our convergence rates to hold is at least $1-4p^{-1}-2p\exp(-\frac{n-1}{16})-c_1\exp(-c_2n)$. With our sample size condition in \eqref{sample_complexity_1_norm}, the above probability has an order of $1-\cO(p^{-1})$.
\end{remark}

\begin{remark}
Note that while the choice of the tuning parameters $c$ and $\lambda$ in \eqref{parameter_setting} guarantees the optimal rates of convergence in both the estimation error and misclassification rate, in practice we recommend to fine-tune these parameters to achieve more appealing performance. In our numerical studies below, we use an independent validation set to tune the parameters in our PANDA method as well as the LPD and AdaLDA method for comparison. We also include the results of our PANDA method with the fine-tuned parameters and with parameters set as in \eqref{parameter_setting} for comparison. 
\end{remark}


\section{Numerical Results}\label{numerical_results}
In this section, we thoroughly compare our proposed PANDA method with the LPD method and AdaLDA method through numerical experiments using both simulated and real data. 

\subsection{Simulation}
To make a fair comparison of the three methods' performances, we fine-tune the parameters for each method on a validation dataset independent from the training data, and we provide both the estimation error of $\beta^*$ (in $\ell_2$ norm) and the population risk  (\ref{risk_formula}) of each method. 

\noindent\textbf{Settings:} 
We follow the settings in \cite{tony2019high} to generate  $\Sigma$ and $\beta^*$. 
\begin{enumerate}[(a)]
	\item \textbf{AR(1)}. We let $\Omega_{j,k} = 0.9^{|j-k|}$, $\Sigma = \Omega^{-1}$ and $\beta^*=(2/\sqrt{s}, \cdots, 2/\sqrt{s},0,\cdots,0)^\top$, where the first $s$ entries are non-zero and $\norm{\beta^*}_2=2$.
	\item \textbf{Varying diagonal}. 
	We let the diagonal entries of $\Sigma$ as $\Sigma_{j,j}=11$ for $j=1,2,\cdots,5$, and $\Sigma_{j,j} = 1+ U_j$ for $j=6,7,\cdots,p$, where $U_i$'s are independently drawn from the uniform distribution $U(0,1)$, and we let the off-diagonal entries  be $\Sigma_{j,k} = 0.9^{|j-k|}$. 
	We let $\beta^* = (1/\sqrt{s},\cdots, 1/\sqrt{s}, 0,\cdots,0)^\top$, where only the first $s$ entries are non-zero and $\norm{\beta^*}_2=1$.
	\item \textbf{Erd{\"o}s-R{\'e}nyi random graph}. We let $\tilde{\Omega}_{j,k} = u_{j,k}v_{j,k}$, where $v_{j,k}$'s are i.i.d. Bernoulli random variables with success rate $0.2$, and $u_{j,k}$'s are i.i.d. uniform random variables  over $[0.5,1] \bigcup[-1,-0.5]$, and $v_{j,k}$'s and $u_{j,k}$'s are  independent. Then we let $\tilde{\Omega}_s = (\tilde{\Omega}+\tilde{\Omega}^\top)/2$ and $\Omega_0 = \tilde{\Omega}_s+\left[\max(-\lambda_{\min}(\tilde{\Omega}_s), 0)+ 0.05\right]I_p$. Let $D_0$ be a diagonal matrix with diagonal elements same as $\Omega_0$'s. We  let $\Omega = D_0^{-1/2}\Omega_0 D_0^{-1/2}$ and $\Sigma = \Omega^{-1}$, and let  $\beta^* = (1/\sqrt{s},\cdots,1/\sqrt{s},0,\cdots,0)^\top$ where only the first $s$ entries are non-zero and $\norm{\beta^*}_2=1$.
	\item \textbf{Block sparse model}. We first construct a matrix $B$ of size $p\times p$ as follows. For $1\leq j\leq p/2$ and $j<k\leq p$,  we let $B_{j,k} = B_{k,j} = 10 b_{i,j}$, where $b_{j,k}$'s are i.i.d. Bernoulli variables with success rate 0.5. For $p/2<j<k\leq p$, we let $B_{j,k} = B_{k,j} = 10$. For the diagonal elements, we let $B_{j,j} = 1$ for $1\leq j\leq p$. Then we let $w = \max(-\lambda_{\min}(B),0) + 0.05$ and let $\Omega = (B + wI_p)/(1+w)$ and $\Sigma = \Omega^{-1}$. We let $\beta^* = (\frac{1}{2\sqrt{s}},\cdots,\frac{1}{2\sqrt{s}},0,\cdots,0)^\top$, where only the first $s$ entries are non-zero and $\norm{\beta^*}_2=1/2$.
	\item \textbf{Approximately sparse setting}. We let $\Sigma_{j,k} = 0.9^{|j-k|}$ and $\beta^*_j = 0.75^j$, 
	which are approximately sparse. Note that $\norm{\beta^*}_2\approx 3$ when $p$ is large.
\end{enumerate}

\begin{figure*}[htb!]
	\centering
	\subfigure[Varying diagonal]{
		\includegraphics[width=0.98\textwidth]{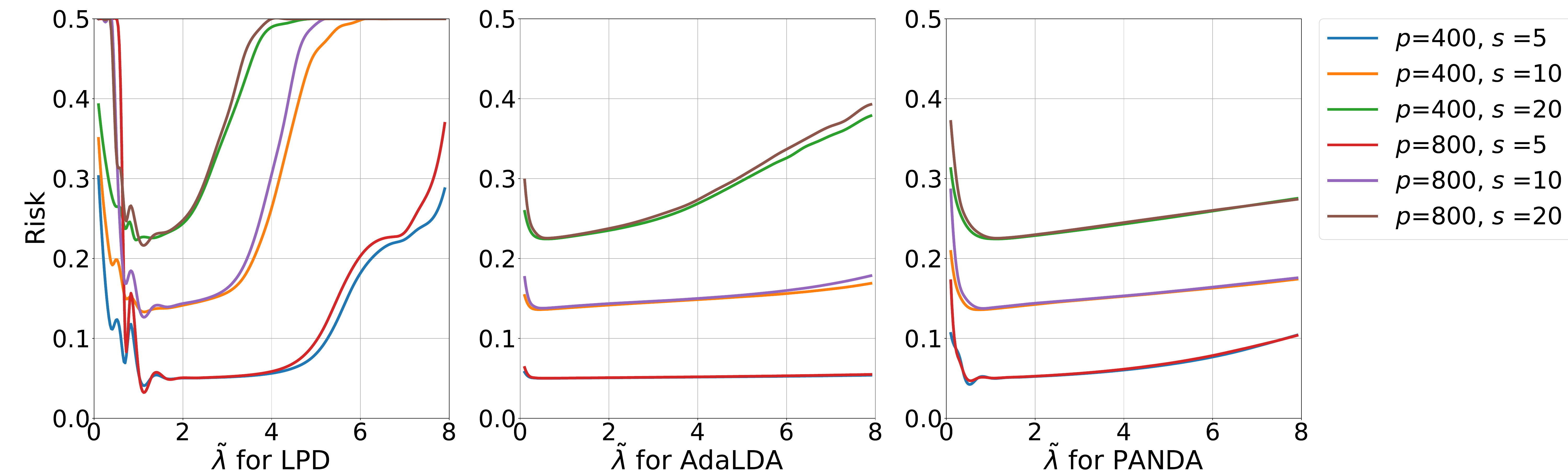}}\\
	\subfigure[Approximately Sparse]{
		\includegraphics[width=0.98\textwidth]{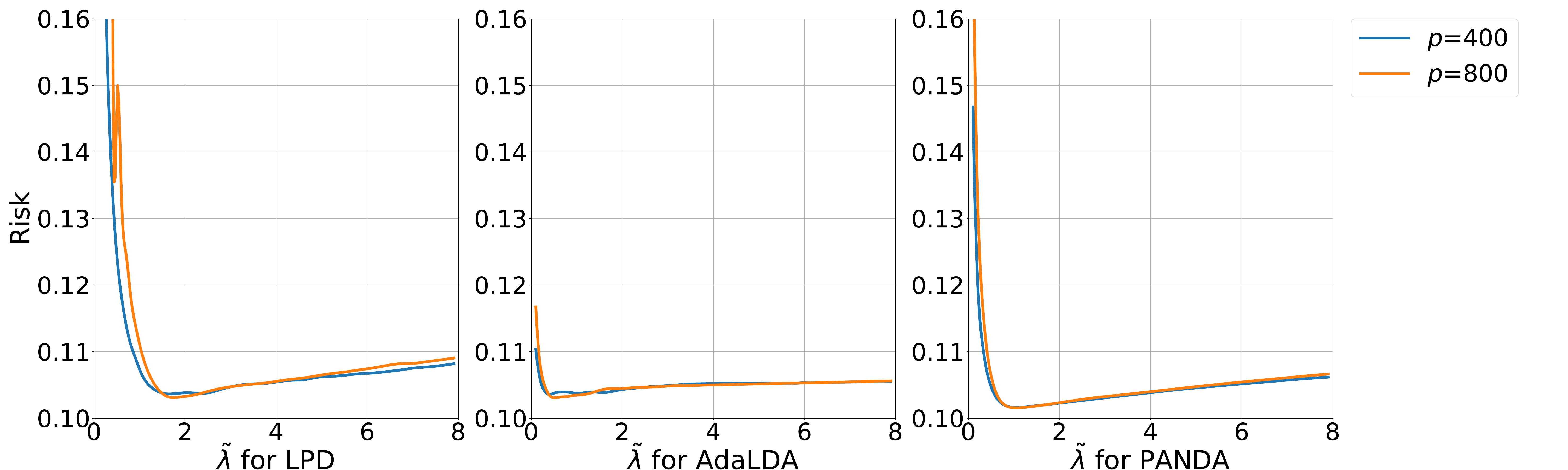}}
	\vspace{-12pt}
	
	\caption{\it The misclassification rate v.s. the values of the parameter $\tilde{\lambda}$ in LPD (left), AdaLDA (middle) and PANDA (right). Results are averaged over 100 replicates.}
	\label{tuning_mis_rate}
	\vspace{-12pt}
\end{figure*}

\begin{figure*}[htb!]
	\centering
	\subfigure[Varying diagonal]{
		\includegraphics[width=0.98\textwidth]{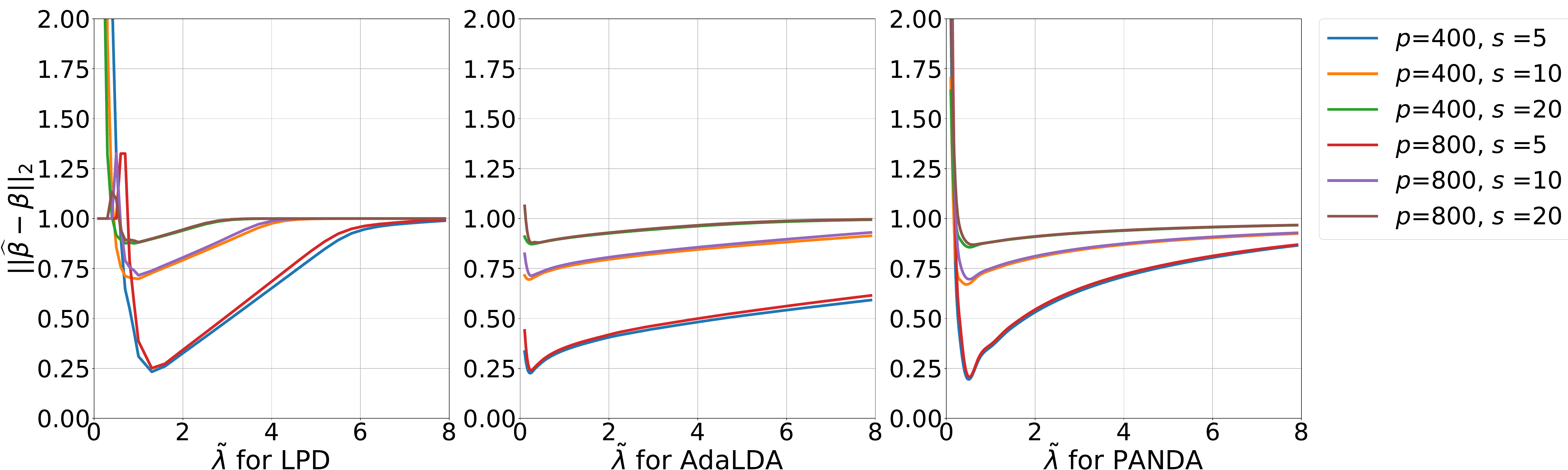}}
	\subfigure[Approximately Sparse]{
		\includegraphics[width=0.98\textwidth]{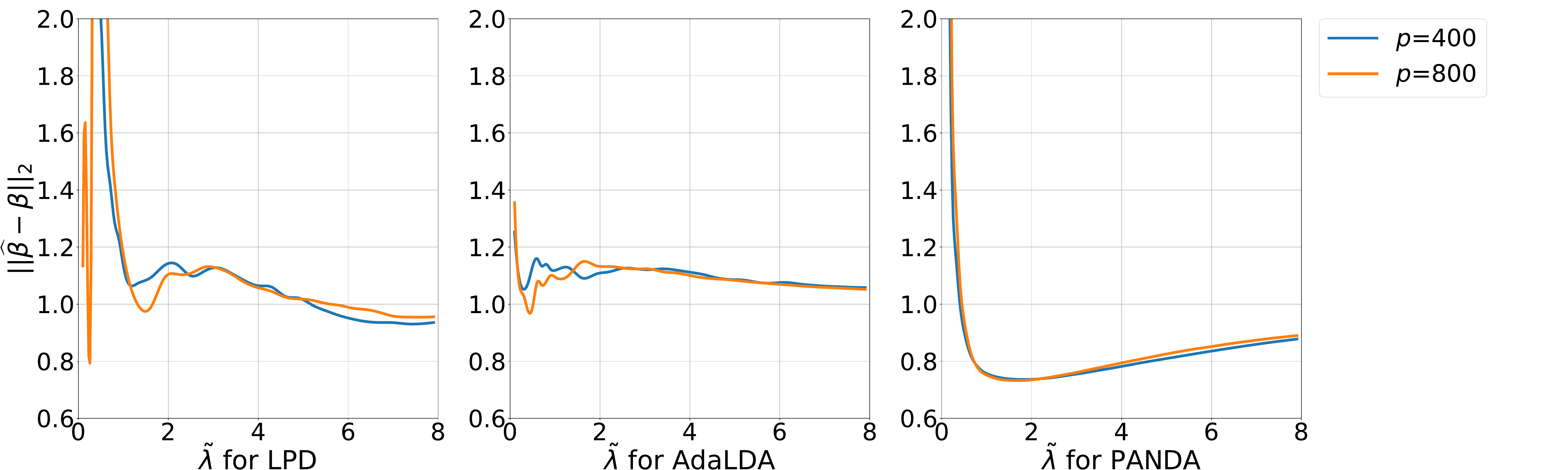}}
	\caption{\it $\ell_2$ estimation error v.s. values of tuning parameter $\tilde{\lambda}$ in LPD (left), AdaLDA (middle) and PANDA (right).
		 Results are averaged over 100 replicates.}
	\label{tuning_error_1}
\end{figure*}

\begin{figure*}[htb!]
	\centering
	\subfigure[Varying diagonal]{
		\includegraphics[width=0.98\textwidth]{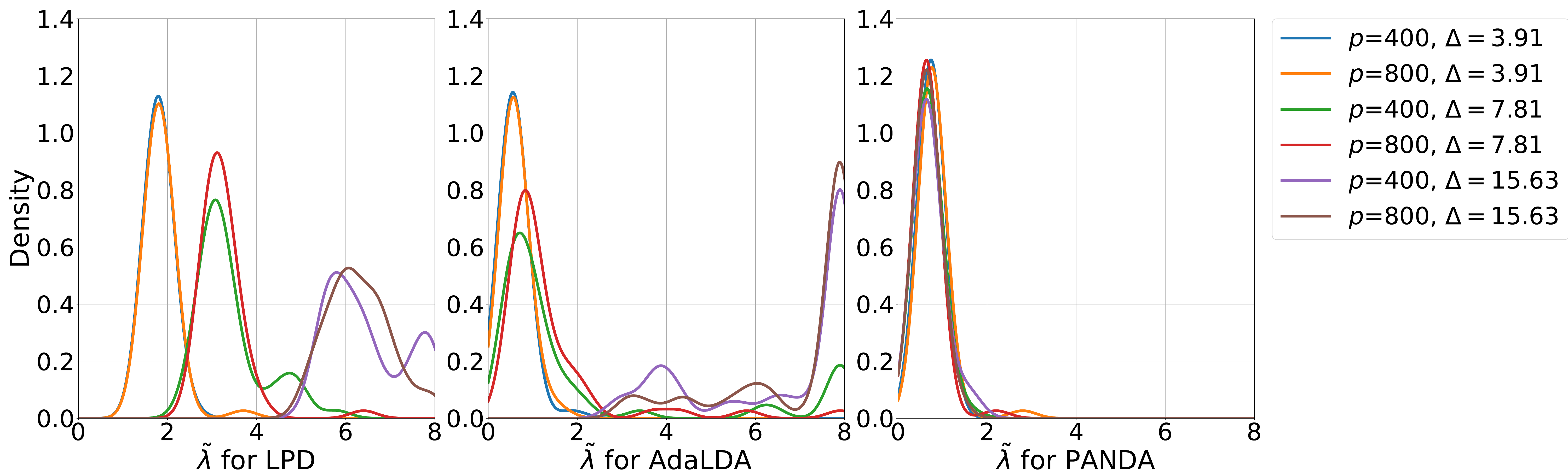}}
	\subfigure[Approximately Sparse]{
		\includegraphics[width=0.98\textwidth]{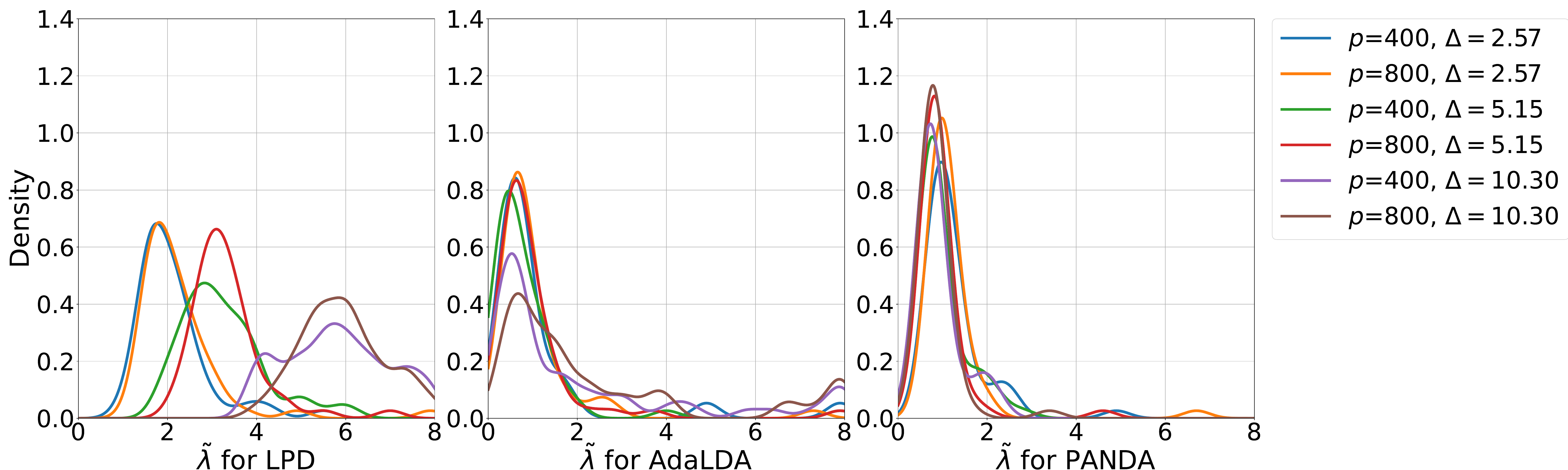}}
	\caption{\it The distribution of the empirically optimal tuning parameter $\tilde{\lambda}$ for LPD (left), AdaLDA (middle) and PANDA (right) over 100 replicates, approximated with kernel smoothing. The optimal choice of the parameter $\tilde{\lambda}$ in our PANDA method relies less on the population.}
	\label{optimal_tuning}
\end{figure*}

\noindent\textbf{Parameter Tuning:}  
While both the AdaLDA method and PANDA method achieve guaranteed theoretical properties with specific tuning parameters, we observe in our experiments that tuning these parameters via a validation set yields better empirical results. In our experiments, under each setting, we randomly sample a validation dataset 
with $n=200$ data points from each class. 
Motivated by the choice of $\lambda$ in (\ref{parameter_setting}), we let $\lambda = \tilde{\lambda}\cdot\sqrt{\log p/n}$, and we tune the parameter $\tilde{\lambda}$, as equivalent to tuning $\lambda$. For a fair comparison, for all the three methods (LPD,  AdaLDA,  and PANDA) we tune $\tilde{\lambda}$  by a grid search over a range from 0.1 to 8.0, with a grid size 0.1.
Figures \ref{tuning_mis_rate} and \ref{tuning_error_1} show the results of the misclassification risks and the estimation errors $\|\hat\beta-\beta^*\|_2$
versus the $\tilde{\lambda}$ value in the three methods, averaged over 100 replicates under each setting of different $p$ and $s$. 
For the parameter $c$ in the PANDA method, we observe that the results are insensitive to the value of $c$ as long as $c$ is not too small, see Table \ref{tab: c insensitive} for the result of the misclassification rate with different choices of $c$  under the AR(1) model as an example. Therefore, we set $c=20$ for all settings. 

\begin{table}[htp!]
	\begin{center}
		\footnotesize
		
		\caption{Misclassification rate of the PANDA method under the AR(1) model with $n=200$, $p=400$, $s=5$ and different $c$,  averaged over 100 replicates. The standard deviations are provided in brackets.}\label{tab: c insensitive}
		
		\begin{tabular}{c|ccc}
			\hline $c$&1e-3&1e-2&0.1\\
			\hline
			Misclassification rate& 
   0.3729 (0.1489)&0.2155 (0.0035)&0.2106 (0.0050)\\
			\hline
			$c$& 1&10&100\\
			\hline
			Misclassification rate&0.2044 (0.0049) &0.2036 (0.0053)&0.2035 (0.0054)\\
			\hline
		\end{tabular}
	\end{center}
	\end{table}

\noindent\textbf{Tuning Sensitivity:} 
We  thoroughly investigate the sensitivity of the tuning parameters under different settings. Since the choice of $\lambda$ in the LPD method relies on the unknown population quantity $\Delta$,  so does the optimal value of $\lambda$ (or $\tilde{\lambda}$, equivalently) in practice.
  We consider following settings to see  how the population distribution, especially the scale of $\Delta$, changes the empirically optimal tuning parameters of  the LPD, AdaLDA and PANDA methods.  
For the varying diagonal model, we set $p=400, 800$, $s=5$, and $\beta^* = \eta\cdot (1/\sqrt{s},\cdots, 1/\sqrt{s}, 0,\cdots,0)^\top$ for $\eta=1,2,4$, where the first $s$ entries are non-zero.
For the approximately sparse $\beta$ model, we set $p=400$, $800$, and $\beta^*_j = \eta\cdot0.75^j$ for $\eta=1,2,4$. 

During the tuning process, we observe that the empirically optimal tuning parameter $\tilde{\lambda}$ for the PANDA method is less sensitive to the change of unknown population quantities among different settings, in comparison with the LPD method and  AdaLDA method. 
In particular, Figure~\ref{optimal_tuning} shows the distribution of the empirically optimal tuning parameter over 100 replicates under each setting
 as specified above. The results show that for the PANDA method, the optimal tuning parameter is always close to $1$, and does not change much across the different settings. 

 \noindent\textbf{Parameter Estimation:} 
Table \ref{table:Estimation} summarizes the estimation error of $\beta^*$, $\|\widehat{\beta}-\beta^*\|_2$,  averaged over 100 random replicates under each setting. It is seen that our proposed PANDA method achieves equal or better performance compared with the LPD and AdaLDA methods in most settings. 

\noindent\textbf{Risk Evaluation:} 
Table \ref{table:Risk} summarizes the misclassification rate under each setting averaged over 100 random replicates. It is seen that our proposed PANDA method achieves similar or better performances than the LPD method and AdaLDA method in most settings. 

\noindent\textbf{Running Time:} Table \ref{VD time} summarizes the running time of our PANDA method and the AdaLDA method under the Varying Diagonal model on a regular computer (Intel Core i5, 2.3GHz). For both methods we use Gurobi, a commercial software that provides state-of-the-art solver for linear programming and second order cone programming, to solve the optimization problems. As can be seen, our PANDA method requires less running time than the AdaLDA method.

\begin{center}
	\begin{table}[htbp!]
			\footnotesize
		\captionsetup{width=1.0\textwidth}
		\caption{\it The $\ell_2$ estimation errors under each setting, averaged over 100 replicates. The standard deviations are given in brackets. The lower value at the significance level $0.05$ between the AdaLDA and the PANDA method are marked in bold.}\label{table:Estimation}
		\begin{tabular}[ht]{p{2.3cm}<{\centering}|cccccccccc}
		\hline \multicolumn{1}{c|}{\textbf{Model}} & \multicolumn{7}{c}{\textbf{Specification}} \\

%
%
%
		\hline
		\multirow{7}*{\makecell{AR(1)\\$\norm{\beta^*}=2$}}&$(s,p)$&$(5, 400)$ & $(10, 400)$ & $(20, 400)$&$(5, 800)$ & $(10, 800)$ & $(20, 800)$ \\
		\cline{2-8}
		&\multirow{2}*{\makecell{LPD}}&1.8875 &1.9607&1.9846&1.8960&1.9669&1.9868\\
		&&(0.0494)&(0.0313)&(0.0101)&(0.0416)&(0.0199)&(0.0094)\\
		&\multirow{2}*{\makecell{AdaLDA}}&1.8854&1.9545 &1.9821&1.8952&1.9593&1.9850\\
		&&(0.0495)&(0.0200)&(0.0098)&(0.0412)&(0.0184)&(0.0084)\\
		&\multirow{2}*{\makecell{PANDA}}&\textbf{1.8673}&\textbf{1.9521}&1.9814&\textbf{1.8856}&1.9571& \textbf{1.9830}\\
		&&(0.0542)&(0.0229)&(0.0112)&(0.0460)&(0.0190)&(0.0104)\\
		
		
		
				\hline
		\multirow{7}*{\makecell{Varying \\Diagonal\\$\norm{\beta^*}=1$}}&$(s,p)$&$(5, 400)$ & $(10, 400)$ & $(20, 400)$&$(5, 800)$ & $(10, 800)$ & $(20, 800)$ \\
		\cline{2-8}
		&\multirow{2}*{LPD}&0.3135 &0.7273&0.8841&0.3158&0.7346&0.8949\\
		&&(0.1088)&(0.0488)&(0.0178)&(0.1128)&(0.0393)&(0.0190)\\
		&\multirow{2}*{AdaLDA} &\textbf{0.2753}  &0.7198 & 0.8837&\textbf{0.2942}&0.7371&0.8935\\
		&&(0.0712)&(0.0387)&(0.0172)&(0.0764)&(0.0374)&(0.0146)\\
		&\multirow{2}*{PANDA} &0.3113& 0.7177&\textbf{0.8797}&0.3197&\textbf{0.7305} &\textbf{0.8901}\\
		&&(0.1110)&(0.0478)&(0.0171)&(0.1166)&(0.0381)&(0.0176)\\
		
		

		\hline
		\multirow{7}*{\makecell{Erd\"{o}s-R\'{e}nyi \\ Random Graph\\
  $\norm{\beta^*}=1$}}&$(s,p)$&$(5, 400)$ & $(10, 400)$ & $(20, 400)$&$(5, 800)$ & $(10, 800)$ & $(20, 800)$& \\
		\cline{2-8}
		& \multirow{2}*{\makecell{LPD}} &0.5715 & 0.7071 & 1.0416 & 0.5933 & 0.7677 & 0.9344\\
		& &(0.1108)&(0.0965)&(0.1608)&(0.1168)&(0.0855)&(0.0867)\\
		&\multirow{2}*{\makecell{AdaLDA}}&0.5688 & 0.6895 & 1.0055 & 0.5949 & 0.7642 & 0.9308 \\
		&&(0.1136)&(0.0761)&(0.0637)&(0.0980)&(0.0914)&(0.1126)\\
		&\multirow{2}*{\makecell{PANDA}} &\textbf{0.5366} & 0.7078 & \textbf{0.9477} & \textbf{0.5753} & \textbf{0.7326} & 0.9114\\
		&&(0.1162)&(0.2120)&(0.0895)&(0.0966)&(0.1054)&(0.2358)\\
		
		\hline
		
		\multirow{7}*{\makecell{Block Sparse\\$\norm{\beta^*}=1$}}&$(s,p)$&$(5, 400)$ & $(10, 400)$ & $(20, 400)$&$(5, 800)$ & $(10, 800)$ & $(20, 800)$ \\
		\cline{2-8}
		&\multirow{2}*{\makecell{LPD}}&0.5066 & 0.5636& 0.6571 & 0.4653 & 0.5490 & 0.5475\\
		&&(0.1184)&(0.1400)&(0.2355)&(0.0908)&(0.0653)&(0.0987)\\
		&\multirow{2}*{\makecell{AdaLDA}}&0.5145 & 0.5480 & 0.5790 & 0.4798 & 0.5391 & \textbf{0.5036} \\
		&&(0.0321)&(0.0082)&(0.0110)&(0.0235)&(0.0143)&(0.0044)\\
		&\multirow{2}*{\makecell{PANDA}}&\textbf{0.4332} & \textbf{0.4986} &\textbf{0.5409} & 0.4789& \textbf{0.5229} & 0.5425 \\
		&&(0.0511)&(0.0272)&(0.0278)&(0.1241)&(0.0665)&(0.1058)	\\
		
		
		\hline		
		\multirow{7}*{\makecell{Approximately\\ Sparse\\$\norm{\beta^*}\approx 3$}}&$p$&$400$ & $800$ 
  & $1200$& \\
		\cline{2-8}
		&\multirow{2}*{\makecell{LPD}}&1.0152 & 0.9900
  & 0.9750 \\
		&&(0.2968)&(0.2897)
  &(0.3112)\\
		&\multirow{2}*{\makecell{AdaLDA}} &1.0117 & 1.0273 
  & 1.0013\\
		&&(0.2877)&(0.2998) 
  &(0.3192)\\
		&\multirow{2}*{\makecell{PANDA}}&\textbf{0.8205} & \textbf{0.8547} 
  & \textbf{0.8514}\\
		&&(0.2328)&(0.2701)
  &(0.2380)\\
		\hline
\end{tabular}
	\end{table}
\end{center}

\begin{center}

	\begin{table}[htbp!]
			\footnotesize
			\captionsetup{width=1.0\textwidth}
			\caption{\it The misclassification rate under each setting averaged over 100 replicates. The standard deviations are given in brackets. The lower value at the significance level $0.05$ between the AdaLDA and the PANDA method are marked in bold.}\label{table:Risk}
			\begin{tabular}{p{2.5cm}<{\centering}|cccccccccc}
					\hline \multicolumn{1}{c|}{\textbf{Model}} & \multicolumn{7}{c}{\textbf{Specification}} \\
			
%
%
%
		\hline
		\multirow{7}*{AR(1)}&$(s,p)$&$(5, 400)$ & $(10, 400)$ & $(20, 400)$&$(5, 800)$ & $(10, 800)$ & $(20, 800)$ \\
		\cline{2-8}
		&\multirow{2}*{LPD}&0.2086 & 0.2900 & 0.3535 & 0.2112 & 0.2908 & 0.3532\\
		&&(0.0074)&(0.0109) &(0.0099)&(0.0074)&(0.0066)&(0.0080)\\
		&\multirow{2}*{AdaLDA} &0.2082  & 0.2890&0.3522&0.2120 &0.2913&\textbf{0.3525}\\
		&&(0.0068)&(0.0080)&(0.0075)&(0.0088)&(0.0072)&(0.0082)\\
		&\multirow{2}*{PANDA} &\textbf{0.2068} & 0.2886 & 0.3542 & 0.2114 & 0.2910 & 0.3571\\
		&&(0.0069)&(0.0087)&(0.0104)&(0.0084)&( 0.0079)&(0.01206)\\
		
		
		
		\hline
		\multirow{7}*{\makecell{Varying\\ Diagonal}}&$(s,p)$&$(5, 400)$ & $(10, 400)$ & $(20, 400)$&$(5, 800)$ & $(10, 800)$ & $(20, 800)$ \\
		\cline{2-8}
		&\multirow{2}*{LPD}&0.0515&0.1382&0.2269&0.0520&0.1390&0.2289\\
		&&(0.0028)&(0.0054)&(0.0065)&(0.0038)&(0.0056)&(0.0087)\\
		&\multirow{2}*{AdaLDA} &0.0508&0.1376&0.2266&0.0513&0.1386&0.2284\\
		&&(0.0018)&(0.0046)&(0.0063)&(0.0032)&(0.0054)&(0.0081)\\
		&\multirow{2}*{PANDA} &0.0512&0.1374&0.2266&0.0514&0.1384&0.2292\\
		&&(0.0026)&(0.0040)&(0.0064)&(0.0025)&(0.0048)&(0.0088)\\
		
		\hline
		\multirow{7}*{\makecell{Erd\"{o}s-R{\'e}nyi\\ Random Graph}}&$(s,p)$&$(5, 400)$ & $(10, 400)$ & $(20, 400)$&$(5, 800)$ & $(10, 800)$ & $(20, 800)$ \\
		\cline{2-8}
		&\multirow{2}*{\makecell{LPD}}&0.2857 & 0.2424 & 0.1150 & 0.2757 & 0.3256 & 0.3289 \\
		&&(0.0138)&(0.0099)&(0.0054)&(0.0148)&(0.0182)&(0.0145)\\
		&\multirow{2}*{\makecell{AdaLDA}} &0.2849 & 0.2414 & 0.1162 & 0.2758 & 0.3246 & 0.3281 \\
		&&(0.0129)&(0.090)&(0.0058)&(0.0138)&(0.0185)&(0.0152)\\
		&\multirow{2}*{\makecell{PANDA}} &\textbf{0.2823} & 0.2403 & \textbf{0.1114} & \textbf{0.2721} & \textbf{0.3183} & \textbf{0.3209}\\
		&&(0.0117)&(0.0106)&(0.0044)&(0.0129)&(0.0166)&(0.0161)\\
		\hline
		\multirow{7}*{Block Sparse}&$(s,p)$&$(5, 400)$ & $(10, 400)$ & $(20, 400)$&$(5, 800)$ & $(10, 800)$ & $(20, 800)$ \\
		\cline{2-8}
		&\multirow{2}*{LPD}& 0.1643 & 0.0954 & 0.0451 & 0.4184 & 0.1724 & 0.3776\\
		&&(0.0056)&(0.0038)&(0.0028)&(0.0170)&(0.0029)&(0.0077)\\
		&\multirow{2}*{AdaLDA} &0.1745 & 0.1002 & 0.0451 & 0.4378 & 0.1739 & 0.3811 \\
		&&(0.0061)&(0.0009)&(0.0003)&(0.0156)&(0.0007)&(0.0020)\\
		&\multirow{2}*{PANDA} &\textbf{0.1614} & \textbf{0.0938} & \textbf{0.0437} & \textbf{0.4168} & \textbf{0.1706} & \textbf{0.3753} \\
		&&(0.0047)&(0.0018)&(0.0007)&(0.0159)&(0.0026)&(0.0072)\\
		
		\hline
		\multirow{7}*{\makecell{Approximately \\Sparse}}&$p$&400&800& 1200 \\
		\cline{2-8}
		&\multirow{2}*{LPD}& 0.1054 & 0.1047 & 0.1053\\
		&&(0.0046)&(0.0030)&(0.0040)\\
		&\multirow{2}*{AdaLDA} & 0.1042 & 0.1043 & 0.1042\\
		&&(0.0029)&(0.0035)&(0.0038)\\
		&\multirow{2}*{PANDA} &0.1034 & 0.1039 & 0.1040 \\
		&&(0.0033)&(0.0038)&(0.0045)\\
		\hline
		\end{tabular}
	\end{table}
\end{center} 

\begin{table}[htb!]
	\centering
	\footnotesize
	\caption{\it Running time  (in seconds) of the PANDA and AdaLDA methods  under the Varying Diagonal model using Gurobi, over 100 replicates. The standard deviations are given in brackets.}
	\label{VD time}
	\begin{tabular}{ccccccc}
		\toprule
		$(s,p)$&$(5, 400)$ & $(10, 400)$ & $(20, 400)$& \\
		\midrule
		AdaLDA &106.739 (2.653) &107.743 (2.588)&107.017 (2.782)\\
		PANDA &70.202 (4.751) & 71.312 (4.389) & 72.112 (4.965)\\
		\midrule
		$(s,p)$&$(5, 800)$ & $(10, 800)$ & $(20, 800)$\\
		\midrule
		AdaLDA &413.262 (13.209) & 413.876 (12.708) &416.793 (12.383)\\
		PANDA &325.486 (16.372) & 326.125 (16.504)& 333.427 (13.554)\\ 
		\bottomrule
	\end{tabular}
\end{table}

	\noindent\textbf{Variable Selection:} We expect our PANDA method is capable for variable selection, as similar to the LPD and AdaLDA method. Here we report the  performance of the three methods in the accuracy of finding the sparse signal, under the AR(1) and Varying Diagonal model as described above. To be more specific, we compute the average of True Positive and True Negative, together with the Precision and Recall for identifying the non-zero entries in $\beta^*$, after applying a threshold at $0.01$ for entries in $\hat{\beta}$. The results under the two models are summarized in Tables \ref{table:variable selection AR_1} and \ref{table:variable selection VD}, respectively.  We see that PANDA achieves comparable performance with LPD and AdaLDA in the sense of accuracy of variable selection. 

	\begin{center}
	\begin{table}[htb!]
			\footnotesize
   \captionsetup{labelfont={color=black},font={color=black}}
		\captionsetup{width=1.0\textwidth}
		\caption{\it The results on variable selection over 100 replicates under the AR(1) model. The standard deviations are given in brackets. }
	\label{table:variable selection AR_1}
		\revise{\begin{tabular}{p{2.5cm}<{\centering}|cccccccccc}
		\hline \multicolumn{1}{c|}{\textbf{Criteria}} & \multicolumn{7}{c}{\textbf{Specification}} \\
		
%
%
%
		\hline
		\multirow{7}*{True Positive}&$(s,p)$&$(5, 400)$ & $(10, 400)$ & $(20, 400)$&$(5, 800)$ & $(10, 800)$ & $(20, 800)$ \\
		\cline{2-8}
		&\multirow{2}*{LPD}&1.95 & 1.51 & 1.27 & 1.76 & 1.17 & 1.02\\
		&&(0.59)&(0.76) &(0.65)&(0.62)&(0.49)&(0.45)\\
		&\multirow{2}*{AdaLDA}& 1.97& 1.60 & 1.34 &1.77  &1.38  & 1.03\\
		&&(0.56)&(0.70) &(0.65)&(0.55)&(0.56)&(0.33)\\
		&\multirow{2}*{PANDA}& 2.20& 1.75 &1.58  & 1.96 &1.51 &1.29 \\
		&&(0.68)&(0.84) &(0.96)&(0.65)&(0.69)&(0.57)\\
		\hline
		\multirow{7}*{True Negative}&$(s,p)$&$(5, 400)$ & $(10, 400)$ & $(20, 400)$&$(5, 800)$ & $(10, 800)$ & $(20, 800)$ \\
		\cline{2-8}
		&\multirow{2}*{LPD}&386.87 &385.1&375.95  &787.72  & 786.57 &778.23 \\
		&&(10.20)&(11.46) &(7.91)&(8.01)&(6.60)&(4.04)\\
		&\multirow{2}*{AdaLDA}&387.55 &385.31  &376.67  &786.91  & 784.42 &777.92 \\
		&&(6.92)&(7.48) &(4.80)&(9.79)&(7.72)&(4.04)\\
		&\multirow{2}*{PANDA}&386.69&384.04  &375.45 &785.80  &783.44 & 775.16 \\
		&&(8.39)&(10.17) &(7.10)&(11.24)&(11.04)&(7.32)\\
		\hline
		\multirow{7}*{Precision}&$(s,p)$&$(5, 400)$ & $(10, 400)$ & $(20, 400)$&$(5, 800)$ & $(10, 800)$ & $(20, 800)$ \\
		\cline{2-8}
		&\multirow{2}*{LPD}&  0.3741& 0.5668&  0.5955&0.4039  &0.6202  &0.7003 \\
		&&(0.2992)&(0.3666) &(0.3778)&( 0.3207)&(0.3726)&(0.3520)\\
		&\multirow{2}*{AdaLDA}&0.3624 &0.5126  & 0.5465 &0.3880  & 0.4406 &0.6379 \\
		&&( 0.2834)&(0.3432) &(0.3514)&(0.3166)&(0.3388)&(0.3537)\\
		&\multirow{2}*{PANDA}& 0.3713 &.4824  &0.4797  &0.4055  & 0.4511 &0.4411 \\
		&&(0.2858)&(0.3373) &(0.3164)&(0.3266)&(0.3427)&(0.3132)\\
		\hline
		\multirow{7}*{Recall}&$(s,p)$&$(5, 400)$ & $(10, 400)$ & $(20, 400)$&$(5, 800)$ & $(10, 800)$ & $(20, 800)$ \\
		\cline{2-8}
		&\multirow{2}*{LPD}&0.3900 & 0.1510 &0.0635  &  0.3520&0.1170  & 0.0510 \\
		&&(0.1185)&(0.07588) &(0.0324)&(0.1243)&(0.0493)&(0.0224)\\
		&\multirow{2}*{AdaLDA}& 0.3940&  0.1600&0.0670 & 0.3540 &0.1380  &0.0515 \\
		&&(0.1118)&(0.0696) &(0.0327)&(0.1096)&(0.0565)&(0.0166)\\
		&\multirow{2}*{PANDA}&0.4400 &  0.1750&  0.0790&0.3920  &  0.1510& 0.0645\\
		&&(0.1363)&(0.0845) &(0.0478)&(0.1300)&(0.0689)&(0.0287)\\
		\hline
		\end{tabular}}
	\end{table}
\end{center}



\begin{center}
\begin{table}[htb!]
		\footnotesize
  \captionsetup{labelfont={color=black},font={color=black}}
	\caption{\it The results on variable selection over 100 replicates under the Varying Diagonal model. The standard deviations are given in brackets. }
	\label{table:variable selection VD}
 \revise{
	\begin{tabular}{p{2.5cm}<{\centering}|cccccccccc}
	\hline \multicolumn{1}{c|}{\textbf{Criteria}} & \multicolumn{7}{c}{\textbf{Specification}} \\
	
%
%
%
	\hline
	\multirow{7}*{True Positive}&$(s,p)$&$(5, 400)$ & $(10, 400)$ & $(20, 400)$&$(5, 800)$ & $(10, 800)$ & $(20, 800)$ \\
	\cline{2-8}
	&\multirow{2}*{LPD}&5.00 & 5.14 & 4.33 & 5.00 & 5.10 & 3.81\\
	&&(0)&(0.85) &(1.14)&(0)&(0.77)&(1.35)\\
	&\multirow{2}*{AdaLDA}& 5.00& 5.25 & 4.34 &5.00  &5.08  & 3.89\\
	&&(0)&(0.54) &(1.12)&(0)&(0.60)&(1.27)\\
	&\multirow{2}*{PANDA}& 5.00& 5.30 &4.44  & 5.00 &5.20 &4.02 \\
	&&(0)&(0.69) &(1.00)&(0)&(0.64)&(1.31)\\
	\hline
	\multirow{7}*{True Negative}&$(s,p)$&$(5, 400)$ & $(10, 400)$ & $(20, 400)$&$(5, 800)$ & $(10, 800)$ & $(20, 800)$ \\
	\cline{2-8}
	&\multirow{2}*{LPD}&394.94 &389.97&380.00  &794.90  & 790.0&780.0 \\
	&&(0.31)&(0.30) &(0)&(0.48)&(0)&(0)\\
	&\multirow{2}*{AdaLDA}&394.95&389.99 &380.0 &794.94 & 789.99 &780.0 \\
	&&(0.26)&(0.10) &(0)&(0.31)&(0.10)&(0)\\
	&\multirow{2}*{PANDA}&394.86 &390.0  &379.98 &794.95  &790.0& 779.99 \\
	&&(0.75)&(0) &(0.20)&(0.26)&(0)&(0.10)\\
	\hline
	\multirow{7}*{Precision}&$(s,p)$&$(5, 400)$ & $(10, 400)$ & $(20, 400)$&$(5, 800)$ & $(10, 800)$ & $(20, 800)$ \\
	\cline{2-8}
	&\multirow{2}*{LPD}& 0.9910& 0.9970&1  & 0.9863 &0.9939  &1 \\
	&&(0.0461)&(0.0302) &(0)&(0.0635)&(0.0010)&(0)\\
	
	&\multirow{2}*{AdaLDA}&0.9921 & 0.9986 &0.9608  & 0.9910 &  0.9983& 1 \\
	&&(0.0400)&(0.0143) &(0.0028)&(0.0461)&(0.0167)&(0)\\
	
	&\multirow{2}*{PANDA}& 0.9830&1  &0.9975  &0.9921  & 1 &0.9985 \\
	&&(0.0733)&(0) &(0.0251)&( 0.0400)&(0)&(0.0145)\\
	\hline
	\multirow{7}*{Recall}&$(s,p)$&$(5, 400)$ & $(10, 400)$ & $(20, 400)$&$(5, 800)$ & $(10, 800)$ & $(20, 800)$ \\
	\cline{2-8}
	&\multirow{2}*{LPD}&1 &  0.514& 0.2165 &1  & 0.51 &0.1905 \\
	&&(0)&(0.0853) &(0.0569)&(0)&(0.0772)&(0.0673)\\
	
	&\multirow{2}*{AdaLDA}& 1&0.525  &  0.217&  1&  0.508&0.1945 \\
	&&(0)&(0.0539) &(0.0560)&(0)&(0.0598)&(0.0635)\\
	
	&\multirow{2}*{PANDA}& 1&  0.53&  0.222&1  &0.52  &0.2010 \\
	&&(0)&(0.0689) &(0.0499)&(0)&(0.0636)&(0.0655)\\
	\hline
	\end{tabular}
 }
\end{table}
\end{center}

\vspace{3cm}

\subsection{Leukemia data}
We investigate the performance of the PANDA, LPD, and AdaLDA methods on a Leukemia dataset from high-density oligonucleotide microarrays. This dataset was first analyzed by \cite{golub1999molecular}, and it contains 72 samples of two categories: 47 of acute lymphoblastic leukemia (ALL), and 25 of acute myeloid leukemia (AML). Each sample contains the quantitative expression levels of 7129 genes.

\noindent\textbf{Preprocessing:} We follow the preprocessing steps  in \cite{tony2019high}. First, we combine the data from both categories and compute the sample variance of each gene. Then, we drop the genes with sample variance beyond the lower and upper 6-quantiles of the total 7129 genes. 

\noindent\textbf{Result:} 
To provide a fair comparison among the LPD, AdaLDA, and PANDA methods, we tune the parameters using a validation set. After preprocessing the raw data, we randomly split the data into training, validation, and testing sets. Specifically, the training set contains 29 ALL and 15 AML samples, the validation set contains 9 ALL and 5 AML samples, and the testing set contains 9 ALL and 5 AML samples. For the computational efficiency, we only use 2000 genes with the largest absolute values of the two-sample $t$-test in the training set, as suggested by \cite{tony2019high}. We repeat the process 100 times, and provide the three methods' average misclassification rates on the testing set (testing error) and their standard deviations in Table \ref{table: Leukemia}. As can be seen, the PANDA method achieves a lower misclassification rate than both the LPD and AdaLDA methods.
	
\begin{table}[htb!]	
	\centering
\caption{\it The performance of PANDA, AdaLDA and LPD on the Leukemia dataset. The testing errors are averaged over 100 replicates. The standard deviation of the testing errors are given in brackets. The difference between PANDA and the other two methods is significant by pair-wise $t$-test with a $p$-value less than 0.001.}
	\vspace{0.1in}
	\begin{tabular}{p{3cm}<{\centering}|ccc}
		\hline
		& LPD& AdaLDA&PANDA\\
		\hline
\multirow{2}*{Testing Error}&9.28\%& 10.64\% &6.93\%\\
&(6.87\%)&(7.92\%)&(6.74\%)	\\
		\hline
	\end{tabular}
	\label{table: Leukemia}
\end{table}


\section{Extension to multiple-class LDA}\label{extension K-class}
	In this section, we discuss the extension of PANDA method to $K$-class LDA in high dimensions. To be more specific, we consider the following data setting. Suppose we have samples $\left\{X^{(k)}_i: k=1,2,\cdots, K, ~i= 1,2,\cdots, n_k\right\}$ from $K$ classes denoted by $k=1,2,\cdots, K$, such that $X^{(k)}_i$'s are i.i.d. from $N(\mu^{(k)}, \Sigma)$. Also, we suppose that the prior probabilities $\pi_1,\pi_2,\cdots,\pi_K$ for the $K$ classes are known.
	Then the oracle classification rule for future data $Z$ is given by $f(Z) = \argmax_{k} D_k$, where $D_1 = 0$, $D_k = \left(Z - \frac{\mu^{(1)} + \mu^{(k)}}{2}\right)^\top \beta^{(k)} + \log\left(\frac{\pi_k}{\pi_1}\right)$, with $\beta^{(k)} = \Sigma^{-1}(\mu^{(k)} - \mu^{(1)})$. In addition, we define $\Delta_k = \sqrt{\beta^{(k)\top}\Sigma\beta^{(k)}}$. Let $\hat{\mu}^{(k)}$ be the sample mean of data in class $k$, and let $\hat{\Sigma}$ be the pooled sample covariance matrix over the $K$ classes. Then, one can construct the classifier by using the $K$-class PANDA method, which simultaneously estimate $\beta^{(k)}$'s and $\Delta_k$'s via the following optimization problems.
	\begin{align}
		(\hat{\beta}^{(k)},\hat{\tau}^k)\in\mathop{\arg\min}_{\beta,\tau}\quad &\norm{\beta}_1 + c_k \tau^2,\label{K-class PANDA optimization}\\
		\text{subject to} \quad&\norm{\hat{\Sigma}\beta - (\hat{\mu}^{(k)} - \hat{\mu}^{(1)})}_\infty\leq \lambda\hat{\sigma}_{\max}(\tau+1),\quad \sqrt{\beta^\top\hat{\Sigma}\beta}\leq \tau.\notag
	\end{align}
Based on $\hat{\beta}^{(k)}$'s, one can construct the classifier by $\hat{f}(Z) = \mathbb{\arg\max}_k~\hat{D}_k$ with $\hat{D}_1 =0$ and $\hat{D}_k = (Z-\frac{\hat{\mu}^{(1)} + \hat{\mu}^{(k)}}{2})^\top \hat{\beta}^{(k)}$.

Following the similar technical argument as for Theorems \ref{one_norm_two_norm_theorem}, \ref{one_norm_two_norm_theorem_revise} and \ref{risk_bound_theorem}, we can establish the following theoretical properties for $K$-class PANDA method.

\begin{theorem}\label{K-class estimation error theorem}
		Suppose that Assumption
	\ref{Sigma_eigenvalue_assumption} hold, and $\beta^{(k)*}\in \mathbb{B}_q(R)$ for some $q\in[0,1)$ and some  $R>0$ for all $k=2,3,\cdots,K$. Let $\left(\widehat{\beta}^{(k)},\widehat{\tau}_k \right)$ be an optimal solution of \eqref{K-class PANDA optimization}. 
	Given
	\begin{align}\label{K-class parameter_setting}
		c_k = \frac{1}{8\left(\norm{\widehat{\mu}^{(k)} - \widehat{\mu}^{(1)}}_\infty + 4\widehat{\sigma}_{\max}\sqrt{\frac{\log p}{n}}\right)}, 
		\quad \lambda = 20 \sqrt{\frac{\log p}{n}},
	\end{align} 
	for sufficiently large $n$ such that 
	\begin{align}\label{K-class sample_complexity_1_norm}
		n\geq C\cdot a^{-2}\Delta_{k}^{2}\sigma_{\max}^2M^{2+\frac{1}{1-q}}R^{\frac{2}{1-q}}\log p
	\end{align}
	where $C$ is an absolute constant, we have, with probability goes to 1,
	\begin{subequations}
		\begin{align}
			\norm{\hat{\beta}^{(k)}-\beta^{(k)*}}_1 &\leq C_1 \cdot (\Delta_k+1)\left(\sigma_{\max} M\right)^{1-q}R\left(\frac{\log p}{n} \right)^{(1-q)/2},\\
			\norm{\hat{\beta}^{(k)}-\beta^{(k)*}}_2 &\leq C_2\cdot (\Delta_k+1) \left(\sigma_{\max} M\right)^{1-q/2}\sqrt{R}\left(\frac{\log p}{n}\right)^{1/2-q/4},\\
			\frac{|\widehat{\tau}_k^2 -\Delta_k^2|}{\Delta_k^2} &\leq C_3\cdot (1+\Delta_k^{-1})\sigma_{\max}^{1-q/2}M^{3/2-q}R\Big(\frac{\log p}{n}\Big)^{(1-q)/2}, 
		\end{align}
	\end{subequations}
	where $C_1$, $C_2$ and $C_3$ are  positive constants.
\end{theorem}
\begin{theorem}
		 Let $\Delta_{\min} = \min
		 \{(\mu^{(j)} - \mu^{(i)})^\top \Sigma^{-1}(\mu^{(j)} - \mu^{(i)}):~1\leq i<j\leq K\}$. Under the identical conditions as in Theorem \ref{K-class estimation error theorem}, we have, with probability goes to 1,
	\begin{align*}
		\mathcal{R}(\hat{f}) - \mathcal{R}^*\leq 
		C\cdot \exp\left(-\frac{\Delta_{\min}^2}{8}\right) \sigma_{\max}^{-q}M^{3-q}\Delta_{\min}
		R \left(\frac{\log p}{n}\right)^{1-q/2},
	\end{align*}
	where $C$ is an absolute positive constants.
\end{theorem}


\section{Proofs of the Main Results}\label{proof}
In this section, we provide the proof for Theorem \ref{one_norm_two_norm_theorem} in Section \ref{proof of theorem 1} and Theorem \ref{risk_bound_theorem} in Section \ref{proof of theorem 2}. The proofs of lemmas can be found in the supplementary material.

\subsection{Proof of Theorem \ref{one_norm_two_norm_theorem}}\label{proof of theorem 1}
\begin{proof}
	We denote by $\delta = \hat{\beta}-\beta^*$ and $\tau^* = \sqrt{\beta^{*\top} \hat{\Sigma}\beta^*}$. We first derive the upper bound for $\norm{\delta}_1$. Based on this upper bound, we then derive the upper bounds for $\norm{\delta}_2$ and $\widehat{\tau}$.
	
	For ease of presentation, we first define the following events,
	\begingroup
	\allowdisplaybreaks
	\begin{align}
		\cE_{\tau} &= \left\{|\beta^{*\top} (\widehat{\Sigma} - \Sigma)\beta^* |\leq \frac{1}{2}\beta^{*\top} \Sigma\beta^* \right\} = \left\{\frac{1}{2}\Delta^2\leq \tau^{*2}\leq \frac{3}{2}\Delta^2 \right\},\label{event_Sigma}\\
		\cE_{\sigma_{\max}} &= \left\{|\widehat{\sigma}^2_{\max} - \sigma^2_{\max}|\leq \frac{1}{2}\sigma^2_{\max} \right\},\label{event_sigma_max}\\
		\cE_{\mu_d} &= \left\{\norm{\mu_d}_\infty - 2\sqrt{2}\sigma_{\max}\sqrt{\frac{\log p}{n}}\leq \norm{\widehat{\mu}_d}_\infty \leq \norm{\mu_d}_\infty + 2\sqrt{2}\sigma_{\max}\sqrt{\frac{\log p}{n}} \right\},\label{event_mu}\\
		\cE_{1} &= \left\{\norm{(\hat{\Sigma} - \Sigma)\beta^*}_\infty \leq 10\sigma_{\max}\Delta \sqrt{\frac{\log p}{n}}\right\},\label{event_1}\\
		\cE_{2} &= \left\{\norm{\hat{\Sigma}\beta^* - \hat{\mu}}_\infty \leq 20 \hat{\sigma}_{\rm{max}}\sqrt{\frac{\log p}{n}} (\tau^* + 1)\right\}.\label{event_2}
		\end{align}
		\endgroup
	
	Before we proceed, we introduce the following lemma.
	\begin{lemma}\label{events_lemma}
	For any $\beta^*\in \mathbb{R}^p$, we have
			\begin{align*}
			\mathbb{P}\left(\cE_\tau \right)&\geq 1-2\exp\left(-\frac{n-1}{16}\right),\quad
			\mathbb{P}\left(\cE_{\sigma_{\max}} \right)\geq 1-2p\exp\left(-\frac{n-1}{16}\right),\\
			\mathbb{P}\left(\cE_{\mu_d} \right)&\geq 1-2p^{-1},\quad
			\mathbb{P}(\cE_1)\geq 1-2p^{-1}.
			\end{align*}
		Moreover, we have
		$$\cE_2\supseteq \left(\cE_{\tau}\bigcap \cE_{\sigma_{\max}}\bigcap \cE_1\right).$$
	\end{lemma}

\noindent\textbf{Upper bound for $\norm{\delta}_1$.}
We first  provide an upper bound 
for $\delta^\top \widehat{\Sigma}\delta$ in terms of $\norm{\delta}_1$,  which is essential for deriving  an upper bound of $\norm{\delta}_1$. 

\begin{lemma}\label{upper_bound_lemma}
	Suppose that 
	the events $\cE_{\tau}$, $\cE_{\sigma_{\max}}$, $\cE_1$ and $\cE_2$ hold.	
	Then we have 
	\begin{align}\label{upper_bound}
		\delta^{\top}\widehat{\Sigma}\delta
		\leq 2\lambda\sigma_{\max}\norm{\delta}_1\left(3\Delta + 2+ \sqrt{\frac{\norm{\delta}_1}{c}}\right).
	\end{align}
\end{lemma}

	
Our next step is to derive a lower bound for $\delta^\top\widehat{\Sigma}\delta$ in terms of $\norm{\delta}_1$, based on the restricted eigenvalue condition of $\widehat{\Sigma}$ on certain restricted subset of $\mathbb{R}^p$. We first introduce the eigenvalue condition of $\widehat{\Sigma}$ that holds with high probability.

	\begin{lemma}\label{Gaussian design RE condition lemma}
	Suppose that Assumption \ref{Sigma_eigenvalue_assumption} holds, and  $n\geq 2$. There exist absolute positive constants $c_1$ and $c_2$ such that 
	\begin{align}\label{Gaussian design RE condition lemma result}
	\delta^\top \widehat{\Sigma}\delta\geq \frac{1}{32M}\norm{\delta}_2^2 - 81\sigma_{\max}^2\frac{\log p}{n}\norm{\delta}_1^2\quad\text{for~all~} \delta\in\mathbb{R}^p,
	\end{align}
	with probability at least $1-c_1\exp(-c_2n)$.
\end{lemma}

Based on the above result, we derive the restricted eigenvalue condition of $\widehat{\Sigma}$ over a restricted subset. In particular, for $S\subseteq[p]$ and $\beta^*\in\mathbb{R}^p$, we let
\begin{align}\label{cone definition}
\cC_{S,\beta^*} \coloneqq \left\{\delta\in\mathbb{R}^p:~ \norm{\delta_{S^c}}_1\leq 3\norm{\delta_S}_1 + 4\norm{\beta^*_{S^c}}_1 \right\}.
\end{align}
The next lemma shows that $\delta\in \cC_{S,\beta^*}$ for any $S\subseteq [p]$.
\begin{lemma}\label{restricted_subset}
	Suppose that 	Assumption \ref{identifiablity_assumption} and 
	events $\cE_{\mu_d}$, $\cE_{\sigma_{\max}}$, $\cE_1$ and $\cE_{2}$, hold. Let $S\subseteq [p]$. Given	$c$ and $\lambda$ in \eqref{parameter_setting}, we have $\delta\in \cC_{S,\beta^*}$ when $n$ satisfies
	\begin{align*}
	n\geq 100 a^{-2}\sigma_{\max}^2\Delta^2\log p.
	\end{align*}
	
\end{lemma}

Now, we choose a subset $S_\eta$ that
\begin{align}
S_\eta &= \left\{j\in[p]:~ |\beta^*_j| \geq \eta \right\},\label{S_choice}\\
\text{where}\quad\eta &= \sigma_{\max}M\sqrt{\frac{\log p}{n}}.\label{eta_choice}
\end{align}
We further show the upper bounds for $|S_\eta|$ and $\norm {\beta^*_{S_\eta^c}}_1$ in the next lemma.
\begin{lemma}\label{weak sparsity upper bounds lemma}
	When $\beta^*\in \mathbb{B}_q(R)$, we have that
\begin{align}
|S_\eta|&\leq \eta^{-q}R,\label{S_eta bound}\\
\norm{\beta^*_{S_\eta^c}}&\leq \eta^{1-q}R.\label{beta_Sc bound}
\end{align}
\end{lemma}

Note that if $S_\eta$ is empty, we immediately have that
\begin{align*}
\norm{\delta}_1 \leq 4\norm{\beta^*_{S_\eta^c}}\leq  4\eta^{1-q}R = 4(\sigma_{\max}M)^{1-q}R\left(\frac{\log p}{n} \right)^{\frac{1-q}{2}},
\end{align*}
which matches the upper bound in (\ref{1-norm_bound}). 

When $S_\eta$ is non-empty and $\delta\in\cC_{S_\eta,\beta^*}$, we have that
\begin{align}\label{1-norm 2-norm}
\norm{\delta}_1\leq 4\norm{\delta_{S_\eta}}_1 + 4\norm{\beta_{S_\eta^c}^*}_1\leq 4\sqrt{|S_\eta|}\norm{\delta}_2 + 4\norm{\beta_{S_\eta^c}^*}_1.
\end{align}
Plugging the above inequality into \eqref{Gaussian design RE condition lemma result} yields that
\begin{align*}
\delta^\top \widehat{\Sigma}\delta&\geq \frac{1}{512M|S_\eta|}\left(\norm{\delta}_1 - 4\norm{\beta^*_{S_\eta^c}}\right)^2 - 81\sigma_{\max}^2\frac{\log p}{n} \norm{\delta}_1^2\\
&\geq \left(\frac{1}{512M|S_\eta|} - 81\sigma_{\max}^2 \frac{\log p}{n}\right)\norm{\delta}_1^2 - \frac{\norm{\beta^*_{S_\eta^c}}_1}{64M|S_\eta|}\norm{\delta}_1.
\end{align*}
When $n$ satisfies that
\begin{align*}
n\geq C\cdot \sigma_{\max}^2 M |S_\eta| \log p
\end{align*}
for some constant $C$, we have that
\begin{align}\label{general RE}
\delta^\top \widehat{\Sigma}\delta\geq \frac{1}{1024M|S_\eta|}\norm{\delta}_1^2 - \frac{\norm{\beta^*_{S_\eta^c}}_1}{64M|S_\eta|}\norm{\delta}_1.
\end{align}
Combining \eqref{upper_bound} with \eqref{general RE}, we have that
\begin{align*}
\frac{1}{1024M|S_\eta|}\norm{\delta}_1^2 - \frac{\norm{\beta^*_{S_\eta^c}}_1}{64M|S_\eta|}\norm{\delta}_1\leq2 \lambda\sigma_{\max}\norm{\delta}_1\left(3\Delta + 2+ \sqrt{\frac{\norm{\delta}_1}{c}}\right).
\end{align*}
Solving the above inequality with our chosen $c$, $\lambda$ and $\eta$ as in \eqref{parameter_setting} and \eqref{eta_choice}, and using the upper bounds \eqref{S_eta bound} and \eqref{beta_Sc bound}, we have the upper bound for $\norm{\delta}_1$ that
\begin{equation}\label{eqn:d1}
\norm{\delta}_1 \leq C\cdot (\sigma_{\max}M)^{1-q}(\Delta + 1)R\left(\frac{\log p}{n}\right)^{\frac{1-q}{2}}
\end{equation}
for some constant $C$, given $n$ satisfies that
\begin{align*}
n\geq C\cdot a^{-2}\Delta^{2}\sigma_{\max}^2M^{2+\frac{1}{1-q}}R^{\frac{2}{1-q}}\log p
\end{align*}
for some constant $C$.

\noindent\textbf{Upper bound for $\norm{\delta}_2$.} We prove (\ref{2-norm_bound}) based on the previous upper bound for $\norm{\delta}_1$. 
Following Lemma \ref{Gaussian design RE condition lemma}, 
there exist some absolute positive constants $c_1$ and $c_2$ such that, with probability at least $1-c_1\exp(-c_2n)$, we have
\begin{align}\label{RE_condition_revisit}
	\delta^\top\widehat{\Sigma}\delta \geq \frac{1}{32M}\norm{\delta}_2^2 -81\sigma_{\max}^2\frac{\log p}{n}\norm{\delta}_1^2.
	\end{align}
	The above inequality gives an upper bound of $\norm{\delta}_2^2$ in terms of $\delta^\top \widehat{\Sigma}\delta$ and $\norm{\delta}_1$, whereas the latter two terms can be further upper bounded using Lemma \ref{upper_bound_lemma} and  \eqref{eqn:d1}, respectively.
	
	
	To bound $\delta^\top \widehat{\Sigma}\delta$, following Lemma \ref{upper_bound_lemma}, we have that
	\begin{align}\label{upper_bound_revisit}
\delta^{\top}\widehat{\Sigma}\delta
\leq \lambda\widehat{\sigma}_{\max}\norm{\delta}_1\left(
\sqrt{\frac{\norm{\delta}_1}{c}} + 3\Delta
+ 2\right).
	\end{align}
	Note that $\Delta = \sqrt{\mu_d^\top \Sigma^{-1}\mu_d}\geq M^{-1/2}\norm{\mu_d}_\infty$, and thus $\norm{\mu_d}_\infty/\Delta \leq M^{1/2}$.
	Hence (\ref{RE_condition_revisit})
	and (\ref{upper_bound_revisit}) together imply that
	\begin{align}
	\norm{\delta}_2^2&\leq 
	C\cdot \left[M \delta^\top \widehat{\Sigma}\delta +  \sigma_{\max}^2M\frac{\log p}{n}\norm{\delta}_1^2\right]\notag\\
	&\leq C\cdot \left[ \lambda\widehat{\sigma}_{\max} M(\Delta+1)\norm{\delta}_1 + \frac{\lambda\widehat{\sigma}_{\max}}{\sqrt{c}}M\norm{\delta}_1^{3/2}
	+ \sigma_{\max}^2M\frac{\log p }{n}\norm{\delta}_1^2\right]\label{2_norm_upper_bound_split}
	\end{align}
	for some constant $C$. By our choice of $c$ and $\lambda$ in (\ref{parameter_setting}) and the upper bound of $\|\delta\|_1$ in (\ref{1-norm_bound}), 
	when $n$ satisfies that
	\begin{equation*}
	n\geq C\cdot 
	\sigma_{\max}^{2}M^{2+\frac{1}{1-q}}R^{\frac{2}{1-q}}\log p
	\end{equation*}
	for some absolute constant $C$, (\ref{2_norm_upper_bound_split}) reduces to 
	\begin{align*}
	\norm{\delta}_2^2 \leq C\cdot \lambda\widehat{\sigma}_{\max} M(\Delta+1)\norm{\delta}_1
	\leq C\cdot \left(\sigma_{\max}M\right)^{2-q}(\Delta+1)^2R\left(\frac{\log p}{n}\right)^{1-q/2},
	\end{align*}
	which shows (\ref{2-norm_bound}) holds. 
	
	\noindent \textbf{Upper bound of $|\widehat{\tau}^2 -\Delta^2|/\Delta^2$. }
	Note that $|\widehat{\tau}^2 -\Delta^2|\leq |\widehat{\tau}^2 - \tau^{*2}| +|\tau^{*2} - \Delta^2|$. We  upper  bound the two terms on the right-hand side respectively in the next lemma. 

		\begin{lemma}\label{Delta_gaps_lemma}
			Suppose that Assumption 
			 \ref{Sigma_eigenvalue_assumption}, events $\cE_{\tau},\cE_{\sigma_{\max}}, \cE_{{\mu}_d}, \cE_{1}$ and (\ref{1-norm_bound}) hold. When $n$ satisfies \eqref{sample_complexity_1_norm}
			for some absolute constant $C$, 
			we have that
			\begin{subequations}
				\begin{align}
				|\widehat{\tau}^2-\tau^{*2}|&\leq C\cdot \Delta(\Delta+1)\sigma_{\max}^{1-q/2}M^{(3-q)/2}R\left(\frac{\log p}{n} \right)^{(1-q)/2} \label{tau_hat_tau*_gap} \\
				|\tau^{*2} - \Delta^2|&\leq  C\cdot\Delta^2 \sigma_{\max}^{1-q/2}M^{(1-q)/2}\sqrt{R}\left(\frac{\log p}{n}\right)^{(2-q)/4} \label{tau*_Delta_gap}
				\end{align}
			\end{subequations}
			for some absolute constant $C$.
		\end{lemma}
	
	Combining (\ref{tau_hat_tau*_gap}) and (\ref{tau*_Delta_gap}), we obtain that
	\begin{equation*}
	\frac{|\widehat{\tau}^2 -\Delta^2|}{\Delta^2}\leq C\cdot (1+\Delta^{-1})\sigma_{\max}^{1-q/2}M^{\frac{3-q}{2}}R\left(\frac{\log p}{n}\right)^{\frac{1-q}{2}} 
	\end{equation*}
	for some absolute constant $C$, and our claim \eqref{Delta_bound} follows as desired.
\end{proof}

\subsection{Proof of Theorem \ref{one_norm_two_norm_theorem_revise}}\label{proof of theorem 1 revised}
\begin{proof}
	We first introduce the following lemma that gives a different upper bound of $\delta^\top \widehat{\Sigma}\delta$ as in Lemma \ref{upper_bound_lemma}, with the additional condition that $\widehat{\tau} = \sqrt{\widehat{\beta}^\top\widehat{\Sigma}\widehat{\beta}}$.
		\begin{lemma}\label{upper_bound_lemma_revise}
			Suppose that 
			the events $\cE_{\tau}$, $\cE_1$ and $\cE_2$ hold, and $\widehat{\tau} = \sqrt{\widehat{\beta}^\top \widehat{\Sigma}\widehat{\beta}}$.	
			Then we have 
			\begin{align}\label{upper_bound_revise}
			\delta^\top \widehat{\Sigma}\delta\leq &C\cdot \lambda \sigma_{\max}\norm {\delta}_1 \Big\{\lambda\sigma_{\max}\norm {\delta}_1 + \tau^*+1+ \left(20\sigma_{\max}\Delta\sqrt{\frac{\log p}{n}}\norm {\delta}_1\right)^{1/2} \notag\\
			&+ \left(2\norm {\mu_d}_2 \norm{\delta}_2\right)^{1/2} \Big\},
			\end{align}
			for some constant $C$.
		\end{lemma}
	
\noindent\textbf{Upper bound of $\norm {\delta}_2$.} Based on Lemma \ref{Gaussian design RE condition lemma} in the previous part, with probability goes to 1 we have that
\begin{align*}
\delta^\top \widehat{\Sigma}\delta\geq \frac{1}{32M}\norm{\delta}_2^2 - 81\sigma_{\max}^2\frac{\log p}{n}\norm{\delta}_1^2\quad\text{for~all~} \delta\in\mathbb{R}^p.
\end{align*}
When $\delta\in \cC_{S_\eta,\beta^*}$, combining the above equation with \eqref{upper_bound_revise}, and using \eqref{1-norm 2-norm}, we have that
\begin{align*}
\frac{1}{M}\norm {\delta}_2^2 \leq C\cdot \left[\sigma_{\max}^2\frac{\log p}{n}\norm {\beta^*_{S_\eta^c}}_1^2 + \lambda\sigma_{\max}\left(\sqrt{|S_\eta|}\norm {\delta}_2+ \norm {\beta^*_{S^c_\eta}}_1\right)\left(\Delta+1+\sqrt{\norm {\mu_d}_2\norm {\delta}_2} \right)\right]
\end{align*}
for some constant $C$, when $n$ satisfies that
\begin{align*}
n\geq C\cdot \sigma_{\max}^2M^{\frac{2-2q}{2-q}}R^{\frac{2}{2-q}}\log p
\end{align*}
for some constant $C$.
By setting $\eta$ as in \eqref{eta_choice}, and using \eqref{S_eta bound} and \eqref{beta_Sc bound}, we finally obtain
\begin{align}\label{2_norm_bound}
\norm{\delta}_2\leq C\cdot (\sigma_{\max}M)^{1-q/2}(\Delta+1)\sqrt{R}\left(\frac{\log p}{n}\right)^{1/2-q/4}
\end{align}
for some constant $C$.

\noindent\textbf{Upper bound of $|\widehat{\tau}^2 - \Delta^2|/\Delta^2$.}	Note that $|\widehat{\tau}^2 -\Delta^2|\leq |\widehat{\tau}^2 - \tau^{*2}| +|\tau^{*2} - \Delta^2|$. In Lemma \ref{Delta_gaps_lemma}, we have already shown the upper bound for the term $|\tau^{*2} - \Delta^2|$ as \eqref{tau*_Delta_gap}, which we also adopt here. With the additional condition that $\sqrt{\widehat{\beta}^\top\widehat{\Sigma}\widehat{\beta}}= \widehat{\tau}$, the upper bound of the term $|\widehat{\tau}^2 - \tau^{*2}|$ can be tighter than \eqref{tau_hat_tau*_gap}, as shown in the following lemma.

		\begin{lemma}\label{Delta_gaps_lemma_revise}
	Suppose that Assumption 
	\ref{Sigma_eigenvalue_assumption}, events $\cE_{\tau},\cE_{\sigma_{\max}}, \cE_{{\mu}_d}, \cE_{1}$ and (\ref{1-norm_bound}) hold. Also,  suppose that $\sqrt{\widehat{\beta}^\top\widehat{\Sigma}\widehat{\beta}^\top} = \widehat{\tau}$. When $n$ satisfies \eqref{sample_complexity_revise},
	we have that
		\begin{align}
		|\widehat{\tau}^2-\tau^{*2}|&\leq C\cdot \Delta(\Delta+1)\sigma_{\max}^{1-q/2}M^{(3-q)/2}\sqrt{R}\left(\frac{\log p}{n} \right)^{1/2-q/4}. \label{tau_hat_tau*_gap_revise}
		\end{align}
	for some absolute constant $C$.
\end{lemma}
Combining \eqref{tau*_Delta_gap} and \eqref{tau_hat_tau*_gap_revise}, we have that
\begin{align*}
\frac{|\widehat{\tau}^2 - \Delta^2|}{\Delta^2}\leq C\cdot (1+\Delta^{-1})\sigma_{\max}^{1-q/2}M^{(3-q)/2}\sqrt{R}\left(\frac{\log p}{n} \right)^{1/2-q/4}. 
\end{align*}
for some constant $C$.
\end{proof}

\subsection{Proof of Theorem \ref{risk_bound_theorem}}\label{proof of theorem 2}
\begin{proof}
	Let 
	$\widehat{\Delta} = \sqrt{\widehat{\beta}^\top \Sigma\widehat{\beta}}$. The misclassification rate of $\widehat{\beta}$ is 
	\begin{align}\label{estimation_risk}
	\mathcal{R}(\widehat{\beta}) = \frac{1}{2}\Phi\left(-\frac{(\widehat{\mu}_m - \mu^{(0)})^\top\widehat{\beta}}{\widehat{\Delta}}\right)  
	+ \frac{1}{2}\Phi\left(\frac{(\widehat{\mu}_m - \mu^{(1)})^\top \widehat{\beta}}{\widehat{\Delta}}\right),
	\end{align}
	where $\Phi(\cdot)$ is the CDF of the standard Gaussian distribution.	Recall that the optimal risk achieved by Fisher's rule is $\cR^* = \Phi(-\frac{\Delta}{2})$. For the first term on the right-hand side of (\ref{estimation_risk}), 
	its second order Taylor's expansion is 
	\begin{align}\label{first_term}
	\Phi\left(-\frac{(\widehat{\mu}_m - \mu^{(0)})^\top\widehat{\beta}}{\widehat{\Delta}}\right) =&~ 
	\Phi\left(-\frac{\Delta}{2}\right) +\Phi'\left(-\frac{\Delta}{2}\right)\left(\frac{\Delta}{2} - \frac{(\widehat{\mu}_m - \mu^{(0)})^\top\widehat{\beta}}{\widehat{\Delta}} \right) \notag\\
	&+ \frac{\Phi''(t_1)}{2} \left(\frac{\Delta}{2} - \frac{(\widehat{\mu}_m - \mu^{(0)})^\top\widehat{\beta}}{\widehat{\Delta}} \right)^2,
	\end{align}
	where $t_1\in\left(\frac{-\Delta}{2},- \frac{(\widehat{\mu}_m - \mu^{(0)})^\top\widehat{\beta}}{\widehat{\Delta}}\right)$.
	Similarly, for the second term in
	(\ref{estimation_risk}), we have
	\begin{align}\label{second_term}
	\Phi\left(\frac{(\widehat{\mu}_m - \mu^{(1)})^\top \widehat{\beta}}{\widehat{\Delta}}\right) =& ~ \Phi\left(-\frac{\Delta}{2}\right) +\Phi'\left(-\frac{\Delta}{2}\right)\left(\frac{\Delta}{2} + \frac{(\widehat{\mu}_m - \mu^{(1)})^\top\widehat{\beta}}{\widehat{\Delta}} \right) 
	\notag\\
	&+ \frac{\Phi''(t_2)}{2} \left(\frac{\Delta}{2} + \frac{(\widehat{\mu}_m - \mu^{(1)})^\top\widehat{\beta}}{\widehat{\Delta}} \right)^2,
	\end{align}
	where $t_2\in(\frac{-\Delta}{2},\frac{(\widehat{\mu}_m - \mu^{(1)})^\top\widehat{\beta}}{\widehat{\Delta}})$.
	Combining (\ref{first_term}) and (\ref{second_term}), we have
	\begin{align}\label{risk_difference}
	\cR(\widehat{\beta}) - \cR^* =& \Phi'\left(-\frac{\Delta}{2}\right) \left(\frac{\Delta}{2} -\frac{\mu_d^\top\widehat{\beta}}{2\widehat{\Delta}}\right)
	+ \frac{\Phi''(t_1)}{2} \left(\frac{\Delta}{2} - \frac{(\widehat{\mu}_m - \mu^{(0)})^\top\widehat{\beta}}{\widehat{\Delta}} \right)^2\notag\\
	&+ \frac{\Phi''(t_2)}{2} \left(\frac{\Delta}{2} + \frac{(\widehat{\mu}_m - \mu^{(1)})^\top\widehat{\beta}}{\widehat{\Delta}} \right)^2.
	\end{align}
	We now introduce a lemma that upper bounds the first term on the right-hand side of \eqref{risk_difference}.
	
	\begin{lemma}\label{first order convergence lemma}
		Suppose \eqref{2-norm_bound} holds, and $n$ satisfies that
		\begin{align*}
		n\geq C\cdot\sigma_{\max}^2 M^{2+2/(2-q)}R^{2/(2-q)}\log p
		\end{align*}
		for some constant $C$. Then we have
		\begin{align}\label{hat_Delta_concentration}
		\frac{\Delta}{2}-\frac{\mu_d^\top\widehat{\beta}}{2\widehat{\Delta}}
		\leq \frac{M}{2\Delta}\norm{\delta}_2^2,
		\end{align}
	\end{lemma}
	
	 Note that $\Phi'(-\Delta/2) = (2\pi)^{-1/2}\exp(-\Delta^2/8)$. Following Lemma \ref{first order convergence lemma}, we have
	 \begin{align}\label{first_term_bound}
	 \Phi^\prime\left(-\frac{\Delta}{2}\right)\left(\frac{\Delta}{2}-\frac{\mu_d^\top\widehat{\beta}}{2\widehat{\Delta}} \right)\leq \frac{M}{2\sqrt{2\pi}\Delta}\exp\left(-\frac{\Delta^2}{8}\right)\norm{\delta}_2^2.
	 \end{align}
	Now we consider the second-order term in (\ref{risk_difference}). 
	First, using Lemma \ref{first order convergence lemma},
	we have
	\begin{align}\label{second_order_term_decompose}
	\frac{\Delta}{2}-\frac{(\widehat{\mu}_m -\mu^{(0)})^\top \widehat{\beta}}{\widehat{\Delta}}&=\frac{\Delta}{2}-\frac{\mu_d^\top\widehat{\beta}}{2\widehat{\Delta}}+\frac{\widehat{\beta}^\top(\mu_m - \widehat{\mu}_m)}{\widehat{\Delta}}\notag\\
	&\leq \frac{M}{2\Delta}\norm{\delta}_2^2 + \frac{\widehat{\beta}^\top (\mu^{(0)} - \widehat{\mu}^{(0)}) + \widehat{\beta}^\top (\mu^{(1)} - \widehat{\mu}^{(1)})}{2\widehat{\Delta}}.
	\end{align}
	After taking square, the first term on the right-hand side gives $\frac{M}{4\Delta^2}\norm{\delta}_2^4$, which is negligible compared to the first-order term. Hence it suffices to bound the second term on the right-hand side of (\ref{second_order_term_decompose}). For this aim we introduce the next lemma.
	\begin{lemma}\label{second order convergence lemma}Under the identical conditions as in Theorem \ref{one_norm_two_norm_theorem} or \ref{one_norm_two_norm_theorem_revise},
		with probability at least $1-4p^{-1}$ we have
		\begin{align}\label{second order term}
		\left(\frac{\widehat{\beta}^\top (\mu^{(0)} - \widehat{\mu}^{(0)}) + \widehat{\beta}^\top (\mu^{(1)} - \widehat{\mu}^{(1)})}{2\widehat{\Delta}}\right)^2
		\leq C\cdot\sigma_{\max}^{-q}M^{1-q}R\left(\frac{\log p}{n}\right)^{1-q/2}
		\end{align}
		for some constant $C$.
	\end{lemma}

	Since $t_1>-\Delta/2$, we have $|\Phi''(t_1)|\leq C\cdot\Delta\exp\left(-\Delta^2/8\right)$.
	Combining this with \eqref{second order term}, we  bound the second term in (\ref{risk_difference}) by
	\begin{align}\label{second_term_bound}
	\frac{|\Phi''(t_1)|}{2} \left(\frac{\Delta}{2} - \frac{(\widehat{\mu}_m - \mu^{(0)})^\top\widehat{\beta}}{\widehat{\Delta}} \right)^2\leq C\cdot \Delta\exp\left(-\frac{\Delta^2}{8}\right)\sigma_{\max}^{-q}M^{1-q}R\left(\frac{\log p}{n} \right)^{1-q/2}
	\end{align}
	for some constant $C$. Likewise, the third term in \eqref{risk_difference} is also subject to this bound.
		
	Finally, plugging \eqref{first_term_bound} and  (\ref{second_term_bound}) into \eqref{risk_difference}, and using \eqref{2-norm_bound}, we achieve that
	\begin{equation*}
	\mathcal{R}(\widehat{\beta}) - \mathcal{R}(\beta^*)\leq
	C\cdot \exp\left(-\frac{\Delta^2}{8}\right) \sigma_{\max}^{-q}M^{3-q}\Delta
	R \left(\frac{\log p}{n}\right)^{1-q/2}
	\end{equation*}
	for some constant $C$, which completes the proof.
\end{proof}

\section{Conclusion and Discussion}\label{discussion}
In this work, we propose PANDA, a novel one-stage and tuning-insensitive method for high-dimensional linear discriminant analysis. We prove that PANDA achieves the optimal convergence rate in both the estimation error and misclassification rate. Our numerical studies show that PANDA achieves equal or better performance compared with existing methods, and requires less effort in parameter tuning. 





Below, we discuss some related work in the existing literature. Besides \cite{gautier2011high}, there are  other pivotal methods for regression and inverse covariance estimation problems. For examples, \cite{belloni2011square} and \cite{sun2012scaled} propose the scaled Lasso method (also known as square-root Lasso) for sparse linear regression, which enjoys a similar tuning-insensitive property to \cite{gautier2011high}; \cite{belloni2014pivotal} extend the scaled Lasso to nonparametric regression; \cite{liu2015calibrated} extend the scaled Lasso to sparse multivariate regression with inhomogeneous noise; \cite{bunea2013group} extend the scaled Lasso to sparse linear regression with group structures; \cite{sun2013sparse} and \cite{liu2017tiger} extend the scaled Lasso to inverse covariance matrix estimation; \cite{zhao2013sparse} extend \cite{gautier2011high} to inverse covariance matrix estimation for heavy tail elliptical distributions; \cite{belloni2011} and \cite{wang2013l1} show that the sparse quantile regression and LAD Lasso are also pivotal methods, which enjoy similar tuning-insensitive properties, respectively.

\acks{The authors thank the action editor and reviewers for their helpful comments, which led to a substantial improvement of the paper. E. X. Fang is partially supported by NSF DMS-2230795 and DMS-2230797. Y. Mei, Y. Shi, and Q. Xu were partially supported by NSF grant DMS-2015405,  NIH grant 1R21AI157618-01A1, and the National Center for
Advancing Translational Sciences of the National Institutes of Health under Award Number UL1TR002378.
The content is solely the responsibility of the authors and does
not necessarily represent the official views of the National Institutes of Health.
}


\newpage

\appendix
\section{An ADMM Algorithm for Solving (\ref{con:problem})}\label{ADMM_scheme}
\label{app:theorem}

This section discusses the implementation of the ADMM algorithm for solving (\ref{con:problem}). For that purpose,
we first re-write the problem (\ref{con:problem}) as
\begin{align}
(\hat{\beta},\hat{\tau}) \in \mathop{\arg\min}_{\beta,u,v,w,
	\in\mathbb{R}^p,\tau
	\in\mathbb{R}}\quad&\norm{\beta}_1+c\tau^2 \label{equiv_problem}\\
\textrm{subject~to}\quad&\widehat{\Sigma}\beta -\lambda\widehat{\sigma}_{\max}\tau
\bm{1}+u = \widehat{\mu}_d + \lambda\widehat{\sigma}_{\max}\bm{1},\notag\\
&\widehat{\Sigma}\beta +\lambda\widehat{\sigma}_{\max}\tau\bm{1} - v= \widehat{\mu}_d - \lambda\widehat{\sigma}_{\max}\bm{1}, \notag\\
&w-\widehat{\Sigma}^{1/2}\beta = 0, \notag\\
&u\geq 0,~v\geq 0, \notag \\
&\norm{w}_2\leq \tau.\notag
\end{align}
Note that the first three constraints in (\ref{equiv_problem}) are linear and the last three constraints  are conic.

To simplify the notation, we write the first three linear constraints as
\begin{equation*}
A_\beta \beta+A_u u +A_v v + A_w w+ A_\tau \tau = b
\end{equation*}
for some real matrices $A_\beta$, $A_u$, $A_v$, $A_w$, $A_\tau$ and real vector $b$. We can further write the problem as
\begin{align*}
(\hat{\beta},\hat{\tau}) \in \mathop{\arg\min}_{\beta,u,v,w
	\in\mathbb{R}^p,\tau
	\in\mathbb{R}}\quad&\norm{\beta}_1+c\tau^2\\
\text{subject~to~}\quad & A_\beta \beta+A_u u +A_v v + A_w w+ A_\tau \tau = b ,\\
& u,v\in \mathcal{C}_1,\\
&(w,\tau)\in \mathcal{C}_2,
\end{align*}
where
\begin{gather*}
\mathcal{C}_1 = \left\{x\in \mathbb{R}^p: x_j\geq 0,~j\in[p] \right\},\\
\mathcal{C}_2 = \left\{(x,y)\in \mathbb{R}^p\times \mathbb{R}: y\geq\sqrt{\sum_{j=1}^p x_j^2} ~\right\}
\end{gather*}
are two convex cones.

The augmented Lagrangian function with scaled dual variables is
\begin{equation*}
L_\rho(\beta,u,v,w,\tau,s) = \norm{\beta}_1 + c\tau^2 + \frac{\rho}{2}\norm{A_\beta \beta + A_u u+A_v v+A_w w+A_\tau \tau -b + s}_2^2-\frac{\rho}{2}\norm{s}_2^2,
\end{equation*}
where $s$ is the scaled dual variable and $\rho>0$ is the penalty parameter.

Based on the  augmented Lagrangian function above, we can derive the ADMM algorithm described in Algorithm \ref{ADMM} in Section \ref{method}. 

In this appendix we prove the following theorem from
Section~6.2:

\noindent
{\bf Theorem} {\it Let $u,v,w$ be discrete variables such that $v, w$ do
not co-occur with $u$ (i.e., $u\neq0\;\Rightarrow \;v=w=0$ in a given
dataset $\dataset$). Let $N_{v0},N_{w0}$ be the number of data points for
which $v=0, w=0$ respectively, and let $I_{uv},I_{uw}$ be the
respective empirical mutual information values based on the sample
$\dataset$. Then
\[
	N_{v0} \;>\; N_{w0}\;\;\Rightarrow\;\;I_{uv} \;\leq\;I_{uw}
\]
with equality only if $u$ is identically 0.} \hfill\BlackBox

\section{Additional Numerical Results}
In this section, we present additional simulation results as supplement to Section \ref{numerical_results}. In subsection \ref{c insensitive}, we include results of PANDA performance with different choices of tuning parameter $c$. In subsection \ref{different n}, we report the performance of LPD, AdaLDA and PANDA when we vary the sample size $n$. In subsection \ref{AUC}, we present the Area Under the Curve (AUC) of the three methods as another performance metric for LDA.

\subsection{PANDA performance with \texorpdfstring{$c$}{c} and \texorpdfstring{$\lambda$}{lambda} in Theorem 1}\label{c insensitive}

In this subsection, we consider the choice of $c$ and $\lambda$ as in \eqref{parameter_setting} for our PANDA method in our simulations. Tables \ref{table: Estimation c} and \ref{table: Risk c} summarizes the performance of our PANDA method with $c$ and $\lambda$ set as in \eqref{parameter_setting}, versus $c=20$ and $\lambda$ fine-tuned under the AR(1) model, together with the performance of LPD and AdaLDA for reference. From these tables, we can see that with parameter $c$ set as in \eqref{parameter_setting}, the PANDA method may not achieve the most desirable empirical performance, and we thus recommend  cross-validation in practice.


	\begin{table}[ht!]
	\begin{center}
 \captionsetup{labelfont={color=black},font={color=black}}
	\caption{\it The $\ell_2$  estimation errors of $\beta^*$ under the AR(1) model, with $n=200$ and different $(s,p)$, averaged over 100 replicates. The standard deviations are given in brackets. The lower value at the significance level $0.05$ between the AdaLDA and the PANDA method are marked in bold.}\label{table: Estimation c}
	\begin{tabular}{
			c|cccccc}
	\hline
	\multirow{2}*{\textbf{Method}}& \multicolumn{6}{c}{$(s,p)$} \\	
	
%
%
%
	&$(5, 400)$ & $(10, 400)$ & $(20, 400)$&$(5, 800)$ & $(10, 800)$ & $(20, 800)$ \\
	\hline
	\multirow{2}*{\makecell{LPD}}&1.8875 &1.9607&1.9846&1.8960&1.9669&1.9868\\
	&(0.0494)&(0.0313)&(0.0101)&(0.0416)&(0.0199)&(0.0094)\\
	\hline
	\multirow{2}*{\makecell{AdaLDA}}&1.8854&1.9545 &1.9821&1.8952&1.9593&1.9850\\
	&(0.0495)&(0.0200)&(0.0098)&(0.0412)&(0.0184)&(0.0084)\\
	\hline
	\multirow{2}*{\makecell{PANDA\\(with $c=20$)}}&\textbf{1.8673}&\textbf{1.9521}&1.9814&\textbf{1.8856}&1.9571& \textbf{1.9830}\\
	&(0.0542)&(0.0229)&(0.0112)&(0.0460)&(0.0190)&(0.0104)\\
	\hline
	\multirow{2}*{\makecell{PANDA\\(with $c,\lambda$ in  Thm 1)} }& 1.9997&2.0000&2.0000&2.0000&2.0000& 2.0000\\
	&( 0.0019)&(0)&(0)&(0)&(0)&()\\
	\hline
\end{tabular}
\end{center}
\end{table}

\begin{table}[ht]
	\begin{center}
 \captionsetup{labelfont={color=black},font={color=black}}
	\caption{\it The misclassification rate under the AR(1) model with different $s$ and $p$, averaged over 100 replicates. The standard deviations are given in brackets.} 
\label{table: Risk c}
\begin{tabular}{
		c|cccccc}
\hline 
\multirow{2}*{\textbf{Method}}& \multicolumn{6}{c}{$(s,p)$} \\

%
%
%
&$(5, 400)$ & $(10, 400)$ & $(20, 400)$&$(5, 800)$ & $(10, 800)$ & $(20, 800)$ \\
\hline
\multirow{2}*{LPD}&0.2086 & 0.2900 & 0.3535 & 0.2112 & 0.2908 & 0.3532\\
&(0.0074)&(0.0109) &(0.0099)&(0.0074)&(0.0066)&(0.0080)\\
\hline
\multirow{2}*{AdaLDA} &0.2082  & 0.2890&0.3522&0.2120 &0.2913&\textbf{0.3525}\\
&(0.0068)&(0.0080)&(0.0075)&(0.0088)&(0.0072)&(0.0082)\\
\hline
\multirow{2}*{\makecell{PANDA\\ (with $c=20$)}} &\textbf{0.2068} & 0.2886 & 0.3542 & 0.2114 & 0.2910 & 0.3571\\
&(0.0069)&(0.0087)&(0.0104)&(0.0084)&( 0.0079)&(0.01206)\\
\hline
\multirow{2}*{\makecell{PANDA\\ (with $c,\lambda$  in Thm 1)}} &0.2444 & 0.3112& 0.3671& 0.2413& 0.3156 & 0.3749\\
&(0.0162)&( 0.0167)&(0.0115)&(0.0165)&(0.0187)&(0.0192)\\
\hline
\end{tabular}
\end{center}
\end{table}

\subsection{Performance of LPD, AdaLDA and PANDA with different \texorpdfstring{$n$}{n}} \label{different n}
Here we present results on the performance of LPD, AdaLDA and our PANDA method with varying sample size. Tables \ref{table: Estimation n} and \ref{table: Risk n} summarize the $\ell_2$ error of $\beta^*$ estimation and the misclassification rate under the AR(1) model, with $n=100$, $200$ and $400$. As can be seen, for every setting of $n$, the three methods achieve comparable performance. 

\begin{table}[ht!]
\captionsetup{labelfont={color=black},font={color=black}}
	\begin{center}
	\caption{\it The $\ell_2$  estimation errors of $\beta^*$ under the AR(1) model, with different $n$, $s$ and $p$, averaged over 100 replicates. The standard deviations are given in brackets. The lower value at the significance level $0.05$ between the AdaLDA and the PANDA method are marked in bold.}\label{table: Estimation n}
\begin{tabular}{p{1.5cm}<{\centering}|cccccccccc}
	\hline
	$n$& \multicolumn{7}{c}{\textbf{Specification}} \\
	\hline
	\multirow{7}*{\makecell{$n=100$}}
	&$(s,p)$&$(5, 400)$ & $(10, 400)$ & $(20, 400)$&$(5, 800)$ & $(10, 800)$ & $(20, 800)$ \\
	\cline{2-8}
	&\multirow{2}*{\makecell{LPD}}&1.9258 & 1.9640&1.9814&1.9236&1.9695&1.9834\\
	&&(0.0408)&(0.0105)&(0.0109)&(0.0396)&(0.0230)&(0.0077)\\
	&\multirow{2}*{\makecell{AdaLDA}}&1.9324&1.9709 &1.9896&1.9298&1.9641&1.9946\\
	&&(0.0292)&(0.0113)&(0.0135)&(0.0326)&(0.0200)&(0.0118)\\
	&\multirow{2}*{\makecell{PANDA} }&1.9161&1.9571&1.9920&\textbf{1.9112}&1.9734& 1.9944\\
	&&(0.0344)&(0.0292)&(0.0199)&(0.0388)&(0.0303)&(0.0140)\\
	\hline
	\multirow{7}*{\makecell{$n=200$}}
	&$(s,p)$&$(5, 400)$ & $(10, 400)$ & $(20, 400)$&$(5, 800)$ & $(10, 800)$ & $(20, 800)$ \\
	\cline{2-8}
		&\multirow{2}*{\makecell{LPD}}&1.8875 &1.9607&1.9846&1.8960&1.9669&1.9868\\
	&&(0.0494)&(0.0313)&(0.0101)&(0.0416)&(0.0199)&(0.0094)\\
	&\multirow{2}*{\makecell{AdaLDA}}&1.8854&1.9545 &1.9821&1.8952&1.9593&1.9850\\
	&&(0.0495)&(0.0200)&(0.0098)&(0.0412)&(0.0184)&(0.0084)\\
	&\multirow{2}*{\makecell{PANDA}}&\textbf{1.8673}&\textbf{1.9521}&1.9814&\textbf{1.8856}&1.9571& \textbf{1.9830}\\
	&&(0.0542)&(0.0229)&(0.0112)&(0.0460)&(0.0190)&(0.0104)\\
			\hline
	\multirow{7}*{\makecell{$n=400$}}
	&$(s,p)$&$(5, 400)$ & $(10, 400)$ & $(20, 400)$&$(5, 800)$ & $(10, 800)$ & $(20, 800)$ \\
	\cline{2-8}
	&\multirow{2}*{\makecell{LPD}}&1.8265 &1.9456 &1.9801&1.8695&1.9824& 3.9300\\
	&&(0.1903) &(0.0247) &(0.0116)&( 0.0601)&(0.0182)&(0.0086)\\
	&\multirow{2}*{\makecell{AdaLDA}}&1.8498 &1.9399 &1.9749& 1.8711&1.9452&1.9775 \\
	&& (0.0764)&(0.0203)&(0.0106)&(0.0370)&(0.0176)&(0.0087)\\
	&\multirow{2}*{\makecell{PANDA}}&\textbf{1.3936} &\textbf{1.9319} &\textbf{1.9706}&\textbf{1.7353}&\textbf{1.9416}&\textbf{1.9748} \\
	&&(0.3866) &(0.0851)&(0.0221)&(0.3031)&(0.0204)&(0.0109)\\
	\hline
\end{tabular}
\end{center}
\end{table}
\begin{table}[ht!]
	\begin{center}
 \captionsetup{labelfont={color=black},font={color=black}}
	\caption{\it The misclassification rate under the AR(1) model, with different $n$, $s$ and $p$, averaged over 100 replicates. The standard deviations are given in brackets. The lower value at the significance level $0.05$ between the AdaLDA and the PANDA method are marked in bold.}\label{table: Risk n}
	\begin{tabular}{p{1.5cm}<{\centering}|ccccccc}
	\hline \multicolumn{1}{c|}{\textbf{$n$}} & \multicolumn{7}{c}{\textbf{Specification}} \\
	\hline
	\multirow{7}*{$n=100$}&$(s,p)$&$(5, 400)$ & $(10, 400)$ & $(20, 400)$&$(5, 800)$ & $(10, 800)$ & $(20, 800)$ \\
	\cline{2-7}
	&\multirow{2}*{LPD}&0.2241 & 0.3019 & 0.3611 & 0.2339 & 0.3152 & 0.3801\\
	&&(0.0092)&(0.0086) &(0.0126)&(0.0097)&(0.0110)&(0.0256)\\
	&\multirow{2}*{AdaLDA} &0.2166  & 0.2969&0.3714&0.2181 &0.3021&\textbf{0.3738}\\
	&&(0.0053)&(0.0068)&(0.0173)&(0.0064)&(0.0093)&(0.0106)\\
	&\multirow{2}*{PANDA} &0.2170 & 0.3136 & 0.3875 & 0.2212 & 0.3214 & 0.4049\\
	&&(0.0082)&(0.0224)&(0.0152)&(0.0077)&( 0.0093)&(0.0206)\\
	\hline
	\multirow{7}*{$n=200$}&$(s,p)$&$(5, 400)$ & $(10, 400)$ & $(20, 400)$&$(5, 800)$ & $(10, 800)$ & $(20, 800)$ \\
	\cline{2-8}
	&\multirow{2}*{LPD}&0.2086 & 0.2900 & 0.3535 & 0.2112 & 0.2908 & 0.3532\\
	&&(0.0074)&(0.0109) &(0.0099)&(0.0074)&(0.0066)&(0.0080)\\
	&\multirow{2}*{AdaLDA} &0.2082  & 0.2890&0.3522&0.2120 &0.2913&\textbf{0.3525}\\
	&&(0.0068)&(0.0080)&(0.0075)&(0.0088)&(0.0072)&(0.0082)\\
	&\multirow{2}*{PANDA} &\textbf{0.2068} & 0.2886 & 0.3542 & 0.2114 & 0.2910 & 0.3571\\
	&&(0.0069)&(0.0087)&(0.0104)&(0.0084)&( 0.0079)&(0.0121)\\
			\hline
\multirow{7}*{\makecell{$n=400$}}
&$(s,p)$&$(5, 400)$ & $(10, 400)$ & $(20, 400)$&$(5, 800)$ & $(10, 800)$ & $(20, 800)$ \\
\cline{2-8}
&\multirow{2}*{\makecell{LPD}}& 0.2000& 0.2815&0.3466&0.2017&0.2824&0.3468\\
&&(0.0056) &(0.0058) &(0.0043)&(0.0058)&(0.0055)&(0.0044)\\
&\multirow{2}*{\makecell{AdaLDA}}&0.1989 &0.2808 &0.3452&0.2003&0.2818&\textbf0.3466\\
&& (0.0042)&(0.0050)& (0.0043)&(0.0042)&(0.0050)&(0.0051)\\
&\multirow{2}*{\makecell{PANDA}}& \textbf{0.1913}& 0.2803&0.3454&0.2000&0.2814&0.3472\\
&& (0.0067)& (0.0053)&(0.0072)&(0.0055)&(0.0059)&(0.0074)\\
	\hline
\end{tabular}
\end{center}
\end{table}

\subsection{AUC of LPD, AdaLDA and PANDA} \label{AUC}
Area Under the Curve (AUC) is another performance metric for binary classification, which looks at the trade-off between the precision and recall rate. In Table \ref{table:AUC} we report the AUC over the testing data with different $s$ and $p$, averaged over 100 replicates. As can be seen, the three methods also achieve comparable performance in AUC.

\begin{table}[p]
	\begin{center}
 \captionsetup{labelfont={color=black},font={color=black}}
	\caption{\it The AUC over testing data,
		averaged over 100 replicates. The standard deviations are given in brackets.} 
\label{table:AUC}
\begin{tabular}{p{2.5cm}<{\centering}|ccccccc}
\hline \multicolumn{1}{c|}{\textbf{Model}} 
& \multicolumn{7}{c}{\textbf{Specification}} \\
\hline
\multirow{7}*{AR(1)}
&$(s,p)$&$(5, 400)$ & $(10, 400)$ & $(20, 400)$&$(5, 800)$ & $(10, 800)$ & $(20, 800)$ \\
\cline{2-8}
&\multirow{2}*{LPD}&0.8770 &0.7858 &0.7034 & 0.8699 & 0.7828 & 0.7051\\
&&(0.0189)&(0.0251) &(0.0297)&(0.0191)&(0.0234)&(0.0295)\\
&\multirow{2}*{AdaLDA} & 0.8773  & 0.7872&0.7048&0.8698 &0.7815&\textbf{0.7059}\\
&&(0.0188)&(0.0238)&(0.0270)&(0.0205)&(0.0228)&(0.0298)\\
&\multirow{2}*{\makecell{PANDA}} & 0.8784 &0.7878 &  0.7028 &0.8700&0.7816 & 0.7001\\
&&(0.0190)&( 0.0245)&(0.0306)&(0.0201)&( 0.0252)&(0.0321)\\

\hline
\multirow{7}*{\makecell{Varying\\ Diagonal}}&$(s,p)$&$(5, 400)$ & $(10, 400)$ & $(20, 400)$&$(5, 800)$ & $(10, 800)$ & $(20, 800)$ \\
\cline{2-8}
&\multirow{2}*{LPD}&0.9898&0.9392&0.8565&0.9899&0.9386&0.8563\\
&&(0.0038)&(0.0125)&(0.0197)&(0.0034)&(0.0109)&(0.0192)\\
&\multirow{2}*{AdaLDA} &0.9899&0.9398&0.8566&0.9900&0.9390&0.8566\\
&&(0.0037)&(0.0119)&(0.0193)&(0.0034)&(0.0106)&(0.0192)\\
&\multirow{2}*{PANDA} &0.9898&0.9401&0.8567&0.9899&0.9390&0.8558\\
&&(.0038)&(0.0117)&(0.0195)&(0.0033)&(0.0108)&(0.0188)\\

\hline
\multirow{7}*{\makecell{Erd\"{o}s-R{\'e}nyi\\ Random Graph}}&$(s,p)$&$(5, 400)$ & $(10, 400)$ & $(20, 400)$&$(5, 800)$ & $(10, 800)$ & $(20, 800)$ \\
\cline{2-8}
&\multirow{2}*{\makecell{LPD}}&0.7826 & 0.8401 & 0.9563 & 0.7992 & 0.7372 & 0.7337 \\
&&(0.0284)&(0.0236)&(0.0101)&(0.0253)&(0.0332)&(0.0257)\\
&\multirow{2}*{\makecell{AdaLDA}} &0.7845 & 0.8415 & 0.9558& 0.7995 & 0.7390 & 0.7353 \\
&&(0.0295)&(0.0241)&(0.0100)&(0.0256)&(0.0325)&(0.0270)\\
&\multirow{2}*{\makecell{PANDA}} &\textbf{0.7867} & 0.8412 & \textbf{0.9589} & \textbf{0.8039} & \textbf{0.7464} & \textbf{0.7439}\\
&&(0.0272)&(0.0236)&(0.0098)&(0.0238)&(0.0316)&(0.0278)\\
\hline
\multirow{7}*{Block Sparse}&$(s,p)$&$(5, 400)$ & $(10, 400)$ & $(20, 400)$&$(5, 800)$ & $(10, 800)$ & $(20, 800)$ \\
\cline{2-7}
&\multirow{2}*{LPD}& 0.9183 & 0.9685 & 0.9920 & 0.6130 & 0.9096 & 0.6688\\
&&(0.0142)&(0.0077)&(0.0034)&(0.0369)&(0.0127)&(0.0280)\\
&\multirow{2}*{AdaLDA} &0.9093 & 0.9660 & 0.9921 & 0.5869 & 0.9082 & 0.6653 \\
&&(0.0156)&(0.0083)&(0.0031)&(0.0331)&(0.0134)&(0.0248)\\
&\multirow{2}*{PANDA} &\textbf{0.9207} & \textbf{0.9696} & \textbf{0.9925} & \textbf{0.6152} & \textbf{0.9113} & \textbf{0.6717} \\
&&(0.0129)&(0.0075)&(0.0031)&(0.0361)&(0.0127)&(0.0291)\\

\hline
\multirow{7}*{\makecell{Approximately \\Sparse}}&$p$&400&800& 1200 \\
\cline{2-8}
&\multirow{2}*{LPD}& 0.9626 & 0.9621 & 0.9624\\
&&(0.0091)&(0.0086)&(0.0075)\\
&\multirow{2}*{AdaLDA} & 0.9625 & 0.9626 & 0.9627\\
&&(0.0098)&(0.0081)&(0.0080)\\
&\multirow{2}*{PANDA} &0.9628 & 0.9621 & 0.9634 \\
&&(0.0098)&(0.0088)&(0.0082)\\
\hline
\end{tabular}
\end{center}
\end{table}

\section{Proofs} \label{all_lemma_proof}

This section provides the detailed proofs to the lemmas in the main body of the paper, and is split into eight subsections, one subsection for the proof of each lemma.

\subsection{Proof of Lemma \ref{events_lemma}}\label{events_lemma_proof}

\begin{proof} There are four main statements in Lemma \ref{events_lemma}, and let us prove them one by one. 

\begin{enumerate}[(i)]
		\item It suffices to show that
		\begin{equation}\label{Delta_concentration_eq_proof}
		\mathbb{P}\left(\cE_\tau \right) = \mathbb{P}\left( |\beta^{*\top} (\widehat{\Sigma} - \Sigma)\beta^* |\leq \frac{1}{2}\beta^{*\top} \Sigma\beta^*\right)\geq 1-2e^{-(n-1)/16}.
		\end{equation}
		Let $\{Y_i\}_{i=1}^{2n-2}$ be i.i.d. random vectors following the multivariate normal distribution  $N({\bf 0}, \Sigma)$. Then
		\begin{equation*}
		\hat{\Sigma}\eqd \frac{1}{2n-2}\sum_{i=1}^{2n-2} Y_i Y_i^\top,\quad \textrm{~and~}\quad
		\beta^{*\top} \hat{\Sigma} \beta^* \eqd \frac{1}{2n-2}\sum_{i=1}^{2n-2} (\beta^{*\top} Y_i)^2,
		\end{equation*}
		where $\eqd$ denotes equal in distribution.
		Note that $\{\beta^{*\top} Y_i\}$ are i.i.d Gaussian random variables following distribution $N(0, \beta^{*\top} \Sigma \beta^*)$, thus $\{(\beta^{*\top} Y_i)^2\}$ are i.i.d. sub-exponential random variables, so for any $t\in (0, \beta^{*\top} \Sigma\beta^*)$, we have
		\begin{equation*}
		\mathbb{P}\left(\Big|\frac{1}{2n-2} \sum_i (\beta^{*\top} Y_i)^2 - \beta^{*\top} \Sigma\beta^*\Big|\geq t\right)\leq 2\exp\left\{-\frac{(2n-2)t^2}{8(\beta^{*\top} \Sigma\beta^*)^2}\right\}.
		\end{equation*}
		Relation (\ref{Delta_concentration_eq_proof}) follows directly by taking $t= \frac{1}{2}\beta^{*\top} \hat{\Sigma}\beta^*$, and thus part (i) of Lemma  \ref{events_lemma} holds. 

		\item Now we need to show that
		\begin{equation*}
		\mathbb{P}\left(\cE_{\sigma_{\max}} \right) =
		\mathbb{P}\left(|\widehat{\sigma}^2_{\max} - \sigma^2_{\max}|\leq \frac{1}{2}\sigma^2_{\max} \right)
		\geq 1-2pe^{-(n-1)/16}.
		\end{equation*}
		To prove this, we set $\beta^* = e_j$ for $j\in [p]$ and use (\ref{Delta_concentration_eq_proof}) with a union bound argument to obtain that
		\begin{equation}
		\mathbb{P}\left(|\hat{\Sigma}_{j,j} - \Sigma_{j,j}|\leq \frac{1}{2}\Sigma_{j,j}, ~\forall j\in [p]\right)\geq 1-2pe^{-(n-1)/16},
		\end{equation}
		where the event on the left-hand side implies that $|\hat{\sigma}^2_{\rm{max}} - \sigma^2_{\rm{max}}|\leq \frac{1}{2}\sigma^2_{\rm{max}}$.

		\item Here it suffices to show that
		\begin{align*}
		\mathbb{P}\left(
			\norm{\mu_d}_\infty - 2\sqrt{2}\sigma_{\max}\sqrt{\frac{\log p}{n}}\leq \norm{\widehat{\mu}_d}_\infty \leq \norm{\mu_d}_\infty + 2\sqrt{2}\sigma_{\max}\sqrt{\frac{\log p}{n}}\right)
		\geq 1-2p^{-1}.
		\end{align*}
		Notice that $\hat{\mu}_d\sim N(\mu_d, \frac{2}{n}\Sigma)$. Let $\mu_{d,j}$ and $\widehat{\mu}_{d,j}$ denote the $j$-th coordinate of $\mu_d$ and $\widehat{\mu}_d$, respectively. We have $\hat{\mu}_{d,j} \sim N(\mu_{d,j}, \frac{2}{n}\Sigma_{j,j})$.
		Therefore, for any $j\in[p]$ we have that
		\begin{equation*}
		\mathbb{P}\left(|\hat{\mu}_{d,j} - \mu_{d,j}| > t \right)\leq 2\exp\left\{\frac{-nt^2}{4(\Sigma_{j,j})^2}\right\}\leq 2\exp\left\{-\frac{nt^2}{4\sigma^2_{\rm{max}}}\right\}.
		\end{equation*}
		Taking $t = \sigma_{\rm{max}}\sqrt{\frac{8\log p}{n}}$ and applying the union bound for all $j\in[p]$, we have with probability at least $1-2p^{-1}$ that
		\begin{equation*}
		|\hat{\mu}_{d,j} - \mu_{d,j}|\leq \sigma_{\rm{max}}\sqrt{\frac{8\log p}{n}},~\forall j\in[p],
		\end{equation*}
		which implies that $|\norm{\hat{\mu}_d}_\infty - \norm{\mu_d}_\infty|\leq 2\sqrt{2}\sigma_{\rm{max}}\sqrt{\log p /n}$.
		\item 	The lower bound of $\mathbb{P}(\cE_{1})$ follows an argument in \cite{tony2019high}. Since $\beta^* = \Sigma^{-1}\mu_d$, we have that $\hat{\Sigma}\beta^* - \hat{\mu}_d = (\hat{\Sigma} - \Sigma)\beta^* - (\hat{\mu}_d - \mu_d)$.
		By A.5.1 in the supplement of \cite{tony2019high}, we have that
		\begin{equation}\label{sigma_beta_infinity_norm}
		\mathbb{P}\left(	|e_j ^\top (\hat{\Sigma} - \Sigma)\beta^* |\leq 10\sigma_{\rm{max}}\Delta \sqrt{\frac{\log p}{n}},~\forall j \in [p] \right)\geq 1-2p^{-1},
		\end{equation}
		where the event on the left-hand side is equivalent to event $\cE_1$. Furthermore, recall that  $\Delta^2 = \beta^{*\top} \Sigma\beta^*$.
		Therefore, under events $\cE_{\tau}$ and $\cE_{\sigma_{\max}}$, we have $\Delta\leq \sqrt{2}\tau^*$ and $\sigma_{\max}\leq \sqrt{2}\widehat{\sigma}_{\max}$. These two conditions and event $\cE_1$ together imply $\cE_2$.
	\end{enumerate}
\end{proof}

\subsection{Proof of Lemma \ref{upper_bound_lemma}}\label{upper_bound_lemma_proof}
\begin{proof}
	When $(\beta^*, \tau^*)$ is feasible to (\ref{con:problem}), from the first constraint of (\ref{con:problem}) we have
	\begin{equation}\label{Sigma delta infinity norm}
	\norm{\hat{\Sigma}\delta}_\infty = \norm{\hat{\Sigma}(\hat{\beta} - \beta^*)}_\infty\leq \norm{\hat{\Sigma}\hat{\beta} - \hat{\mu}_d}_\infty + \norm{\hat{\Sigma}\beta^* - \hat{\mu}_d}_\infty \leq \lambda\widehat{\sigma}_{\max}(\hat{\tau} + \tau^*) + 2\lambda\widehat{\sigma}_{\max}.
	\end{equation}
	In addition, due to the optimality of $(\widehat{\beta},\widehat{\tau})$, we have
	\begin{align*}
	\norm{\widehat{\beta}}_1 + c\widehat{\tau}^2 \leq \norm{\beta^*}_1 + c\tau^{*2},
	\end{align*}
	which implies that
	\begin{align*}
	\widehat{\tau}\leq \tau^* + \sqrt{\frac{\norm{\delta}_1}{c}}.
	\end{align*}
	Plugging the above inequality into (\ref{Sigma delta infinity norm}), we obtain that
	\begin{align*}
	\norm{\hat{\Sigma}\delta}_\infty\leq 2\lambda\widehat{\sigma}_{\max}(\tau^*+1) + \lambda\widehat{\sigma}_{\max}\sqrt{\frac{\norm{\delta}_1}{c}}.
	\end{align*}
	Under the events $\cE_{\tau}$ and $\cE_{\sigma_{\max}}$, we have $\tau^*\leq \sqrt{\frac{3}{2}}\Delta$ and $\widehat{\sigma}_{\max}\leq 2\sigma_{\max}$, so we further have that
	\begin{align*}
	\norm{\hat{\Sigma}\delta}_\infty\leq 2\lambda\sigma_{\max}\left(3\Delta + 2+ \sqrt{\frac{\norm{\delta}_1}{c}}\right).
	\end{align*}
	Finally, applying H{\" o}lder's inequality, we obtain that
	\begin{align*}
	\delta^\top \widehat{\Sigma}\delta\leq \norm{\delta}_1\norm{\widehat{\Sigma}\delta}_\infty
	\leq 2
	\lambda\sigma_{\max}\norm{\delta}_1\left(3\Delta + 2+ \sqrt{\frac{\norm{\delta}_1}{c}}\right).	
	\end{align*}
	Thus Lemma \ref{upper_bound_lemma} holds.
\end{proof}

\subsection{Proof of Lemma \ref{Gaussian design RE condition lemma}}\label{Gaussian design RE condition_lemma_proof}
\begin{proof} Lemma \ref{Gaussian design RE condition lemma} is an application of a theorem in \cite{raskutti2010restricted}, which is given by the following lemma.
	\begin{lemma}[Theorem 1 of \cite{raskutti2010restricted}]\label{Gaussian random design RE}
		For any Gaussian random design $Z\in\mathbb{R}^{n\times p}$ with i.i.d. $N(\bm{0},\Sigma)$ raws, there exist absolute positive constants $c_1,c_2$ such that
		\begin{align*}
		\frac{\norm{Z\delta}_2}{\sqrt{n}}\geq \frac{1}{4}\norm{\Sigma^{1/2}\delta}_2 - 9\sigma_{\max}\sqrt{\frac{\log p}{n}}\norm{\delta}_1,~\forall \delta\in\mathbb{R}^p,
		\end{align*}
		with probability at least $1-c_1\exp(-c_2 n)$.
	\end{lemma}
	
	Now we are ready to prove Lemma \ref{Gaussian design RE condition lemma}. Suppose $n\geq 2$. Then the pooled covariance matrix $\widehat{\Sigma}$ is obtained by
	\begin{equation*}
	\widehat{\Sigma} = \frac{1}{2n-2}\left[\sum_{i=1}^n \left(X^{(0)}_i - \widehat{\mu}^{(0)}\right)\left(X^{(0)}_i - \widehat{\mu}^{(0)}\right)^\top +  \sum_{i=1}^n \left(X^{(1)}_i - \widehat{\mu}^{(1)}\right)\left(X^{(1)}_i - \widehat{\mu}^{(1)}\right)^\top \right],
	\end{equation*}
	and   $\widehat{\Sigma}$ has the same distribution as
	\begin{equation*}
	\tilde{\Sigma} = \frac{1}{2n-2}\sum_{i=1}^{2n-2} Z_i Z_i^\top,
	\end{equation*}
	where $Z_j$'s are i.i.d. samples from $N(0, \Sigma)$. Hence $\widehat{\Sigma}$ can be viewed as the sample covariance matrix of a Gaussian random design with $0$ mean.
	
	By Lemma \ref{Gaussian random design RE} (i.e., Theorem 1 of \cite{raskutti2010restricted}), there exist absolute positive constants $c_1$ and $c_2$ such that with probability at least $1-c_1 \exp(-c_2n)$,
	\begin{equation*}
	\norm{\widehat{\Sigma}^{1/2}\delta}_2\geq \frac{1}{4}\norm{\Sigma^{1/2}\delta}_2 - 9\sigma_{\max}\sqrt{\frac{\log p}{2n-2}}\norm{\delta}_1.
	\end{equation*}	
	When $n\geq 2$ and $\lambda_{\min}\geq M^{-1}$, we have
	\begin{align*}
	\norm{\widehat{\Sigma}^{1/2}\delta}_2\geq \frac{1}{4\sqrt{M}}\norm{\delta}_2 - 9\sigma_{\max}\sqrt{\frac{\log p}{n}}\norm{\delta}_1,
	\end{align*}
	and thus
	\begin{align*}
	\delta^\top \widehat{\Sigma}\delta&\geq \left(\frac{1}{4\sqrt{M}}\norm{\delta}_2 - 9\sigma_{\max}\sqrt{\frac{\log p}{n}}\norm{\delta}_1\right)^2\\
	&\geq\frac{1}{32M}\norm{\delta}_2^2 - 81\sigma_{\max}^2\frac{\log p}{n}\norm{\delta}_1^2.
	\end{align*}
	Here the last inequality follows from the fact that 
	\begin{align*}
	(a-b)^2 = \left(\frac{1}{2}a^2-2ab + 2b^2\right) + \frac{1}{2}a^2 - b^2\geq  \frac{1}{2}a^2 - b^2
	\end{align*}
	for any number $a,b\geq 0$. Thus  Lemma \ref{Gaussian design RE condition lemma} holds.
\end{proof}

\subsection{Proof of Lemma \ref{restricted_subset}}\label{restricted_subset_proof}
\begin{proof}
	For any $S\subseteq [p]$, we have
	\begin{align*}
	\norm{\widehat{\beta}}_1 = \norm{\beta^*+\delta}_1\geq\norm{\beta_S^*}_1+\norm{\delta_{S^c}}_1 - \norm{\beta^*_{S^c}}_1 -\norm{\delta_S}_1.
	\end{align*}
	Combining the above inequality with $\norm{\beta^*}_1 \leq \norm{\beta^*_S}_1 + \norm{\beta^*_{S^c}}_1$, we have
	\begin{align}\label{deviation ineq}
	\norm{\widehat{\beta}}_1-\norm{\beta^*}_1\geq \norm{\delta_{S^c}}_1 - \norm{\delta_S}_1 - 2\norm{\beta^*_{S^c}}_1.
	\end{align}
	When $(\beta^*, \tau^*)$ is feasible to \eqref{con:problem}, by optimality we have
	\begin{equation}
	\norm{\hat{\beta}}_1 + c\hat{\tau}^2\leq \norm{\beta^*}_1 + c\tau^{*2}.
	\label{optim}
	\end{equation}
	Combining \eqref{deviation ineq} and \eqref{optim} yields
	\begin{align}
	\norm{\delta_{S^c}}_1 -\norm{\delta_S}_1 -2\norm{\beta^*_{S^c}}_1\leq \norm{\widehat{\beta}}_1-\norm{\beta^*}_1\leq c(\tau^{*2}-\widehat{\tau}^2).\label{tau diff}
		\end{align}
	Since $\tau^{*2} = \beta^{*\top}\widehat{\Sigma}\beta^*$ and $\widehat{\tau}^2\geq \widehat{\beta}^\top\widehat{\Sigma}\widehat{\beta}$,
	it follows that
	\begin{align*}
	c(\tau^{*2} - \widehat{\tau}^2)
	&\leq -2c\delta^\top(\widehat{\Sigma}\beta^*)\\
	&= -2c\delta^\top (\widehat{\Sigma}-\Sigma)\beta^* -2c\delta^\top\mu_d\\
	&\leq 2c\norm{(\widehat{\Sigma}-\Sigma)\beta^*}_\infty\norm{\delta}_1+2c\norm{\mu_d}_\infty \norm{\delta}_1.
	\end{align*}
	Under event $\cE_1$, we have
	\begin{align*}
	\norm{(\widehat{\Sigma}-\Sigma)\beta^*}_\infty\leq 10\sigma_{\max}\Delta\sqrt{\frac{\log p}{n}}.
	\end{align*}
	When $n$ satisfies that
	\begin{align*}
	n\geq 100 a^{-2}\sigma_{\max}^2\Delta^2\log p,
	\end{align*}
	we have
	\begin{align*}
	c(\tau^{*2} - \widehat{\tau}^2)\leq 4c\norm{\mu_d}_\infty\norm{\delta}_1.
	\end{align*}
%
%
%
	By setting $c$ as in (\ref{parameter_setting}), we have that 
	$$\frac{1}{2}\norm{\delta_{S^c}}_1\leq \frac{3}{2} \norm{\delta_S}_1 + 2\norm{\beta^*_{S^c}}_1.$$
	Thus $\norm{\delta_{S^c}}_1\leq 3\norm{\delta_S}_1 + 4\norm{\beta^*}_1$, which completes the proof of Lemma \ref{restricted_subset}.
\end{proof}

\subsection{Proof of Lemma \ref{weak sparsity upper bounds lemma}}\label{weak sparsity upper bounds lemm proof}
\begin{proof}
		From the definitions of $\mathbb{B}_q(R)$ and $S_\eta$, we have that 
		\begin{equation*}
		R\geq\sum_{j}|\beta^*_j|^q\geq \eta^q |S_\eta|, 
		\end{equation*}
		and
		\begin{equation*}
		R\geq\sum_j |\beta^*_j|^q = \sum_j|\beta^*_j|\cdot |\beta^*_j|^{q-1}\geq \eta^{q-1} \norm{\beta^*_{S_\eta^c}}_1.
		\end{equation*}
		Lemma \ref{weak sparsity upper bounds lemma} follows immediately from these two inequalities, and thus holds. 
\end{proof}

\subsection{Proof of Lemma \ref{Delta_gaps_lemma}}\label{Delta_gaps_lemma_proof}
\begin{proof} Let us first prove relation (\ref{tau_hat_tau*_gap}). Under the optimality condition, we have $\norm{\widehat{\beta}}_1 + c\widehat{\tau}^2 \leq \norm{\beta^*}_1 + c\tau^{*2}$, and thus
	\begin{equation}
	\widehat{\tau}^2-\tau^{*2}\leq \frac{1}{c}\norm{\delta}_1\leq C\cdot \Delta(\Delta+1)\sigma_{\max}^{1-q}M^{3/2-q}R\left(\frac{\log p}{n}\right)^{\frac{1-q}{2}} \label{tau_hat upper bound}
	\end{equation}
	for some positive constant $C$. Here the last inequality uses (\ref{1-norm_bound}) and $\norm{\mu_d}_\infty \leq M^{1/2}\Delta$.
	
	Note that the second constraint in (\ref{con:problem}) implies that
	\begin{align*}
	\widehat{\tau}^2 \geq \widehat{\beta}^\top\widehat{\Sigma}\widehat{\beta}
	= (\beta^*+\delta)^\top\widehat{\Sigma}(\beta^*+\delta)
	\geq \tau^{*2} + 2\delta^\top \widehat{\Sigma} \beta^*,
	\end{align*}
	hence
	\begin{align}
	\widehat{\tau}^2 - \tau^{*2}&\geq -2|\delta^\top\widehat{\Sigma}\beta^*|\\
	&\geq -2\left|\delta^\top \left[(\widehat{\Sigma}-\Sigma)\beta^* + \mu_d\right]\right|\\
	&\geq -2\norm{\delta}_2\norm{\mu_d}_2 -2\norm{\delta}_1\norm{(\widehat{\Sigma} - \Sigma)\beta^*}_\infty.\label{tau_hat_lower_bound}
	\end{align}
	Note that under the event $\cE_1$, we have
	\begin{align}\label{E_1 event}
	\norm{(\widehat{\Sigma}-\Sigma)\beta^*}_\infty\leq 10\sigma_{\max}\Delta\sqrt{\frac{\log p}{n}}.
	\end{align}
	 Plugging \eqref{1-norm_bound}, \eqref{2-norm_bound}, \eqref{E_1 event} and $\norm{\mu_d}_2\leq M^{1/2}\Delta$ into \eqref{tau_hat_lower_bound}, we obtain that
	 \begin{align*}
	 \widehat{\tau}^2 - \tau^{*2}\geq -C\cdot \Delta(\Delta+1) \sigma_{\max}^{1-q/2}M^{(3-q)/2} \sqrt{R}\left(\frac{\log p}{n}\right)^{1/2-q/4}.
	 \end{align*}
	 Combining the above equation and (\ref{tau_hat upper bound}) yields (\ref{tau_hat_tau*_gap}).

	Next, let us prove the result (\ref{tau*_Delta_gap}) in Lemma \ref{Delta_gaps_lemma}. Note that the gap between $\tau^{*2}$ and $\Delta^2$ can be written as
	$|\tau^{*2} - \Delta^2| = |\beta^{*\top}(\widehat{\Sigma}- \Sigma)\beta^*|$. To bound this gap, we first apply H{\"o}lder's inequality that
	\begin{align*}
	|\beta^{*\top}(\widehat{\Sigma}- \Sigma)\beta^*|\leq \norm{\beta^*}_1\norm{(\widehat{\Sigma}- \Sigma)\beta^*}_\infty.
	\end{align*}
	Under event $\cE_1$, the term $\norm{(\widehat{\Sigma}- \Sigma)\beta^*}_\infty$ can be again bounded by (\ref{E_1 event}).
	To bound the term $\norm{\beta^*}_1$, we note that
	\begin{align*}
	\norm{\beta^*}_1 = \norm{\beta^*_{S_\eta}}_1 + \norm{\beta^*_{S_\eta^c}}_1\leq \sqrt{|S_\eta|}~\norm{\beta^*}_2 + \norm{\beta^*_{S_\eta^c}}_1\leq \eta^{-q/2}\sqrt{R}M^{1/2}\Delta + \eta^{1-q}R.
	\end{align*}
	The last inequality above uses equations \eqref{S_eta bound} and \eqref{beta_Sc bound}.
	By our choice of $\eta$ in (\ref{eta_choice}), when $n$ satisfies that
	\begin{align*}
	n\geq C\cdot \Delta^2\sigma_{\max}^2MR\log p
	\end{align*}
	for some absolute constant $C$, we have that 
	\begin{align*}
	\norm{\beta^*}_1\leq C\cdot \eta^{-q/2}\sqrt{R}M^{1/2}\Delta\leq C\cdot\Delta^2 \sigma_{\max}^{-q/2}M^{(1-q)/2}\Delta\sqrt{R}\left(\frac{\log p}{n}\right)^{-q/4}.
	\end{align*}
	Hence we have
	\begin{align*}
	|\tau^{*2} - \Delta^2| = |\beta^{*\top}(\widehat{\Sigma}- \Sigma)\beta^*|\leq \norm{\beta^*}_1\norm{(\widehat{\Sigma}- \Sigma)\beta^*}_\infty \leq C\cdot \sigma_{\max}^{1-q/2}M^{(1-q)/2} \sqrt{R}\left(\frac{\log p}{n}\right)^{\frac{2-q}{4}},
	\end{align*}
	and thus \eqref{tau*_Delta_gap} holds. 
\end{proof}

	\subsection{Proof of Lemma \ref{upper_bound_lemma_revise}}\label{upper_bound_lemma_revise_proof}
	\begin{proof}
		When $(\beta^*, \tau^*)$ is feasible to (\ref{con:problem}), from the first constraint of (\ref{con:problem}) we have
		\begin{equation}\label{Sigma delta infinity norm revisit}
		\norm{\hat{\Sigma}\delta}_\infty = \norm{\hat{\Sigma}(\hat{\beta} - \beta^*)}_\infty\leq \norm{\hat{\Sigma}\hat{\beta} - \hat{\mu}_d}_\infty + \norm{\hat{\Sigma}\beta^* - \hat{\mu}_d}_\infty \leq \lambda\widehat{\sigma}_{\max}(\hat{\tau} + \tau^*) + 2\lambda\widehat{\sigma}_{\max}.
		\end{equation}
		When $\widehat{\tau} = \sqrt{\widehat{\beta}^\top \widehat{\Sigma}\widehat{\beta}}$, we have
		\begin{align*}
		\widehat{\tau}^2 &= \widehat{\beta}^\top \widehat{\Sigma}\widehat{\beta}
		= (\beta^* + \delta)^\top \widehat{\Sigma}(\beta^*+\delta)
		= \tau^{*2} + 2\delta^\top \widehat{\Sigma}\beta^* + \delta^\top \widehat{\Sigma}\delta\\
		&= \tau^{*2}+ \delta^\top \widehat{\Sigma}\delta +2\delta^\top (\widehat{\Sigma} - \Sigma)\beta^* + 2\delta^\top \mu_d\\
		&\leq \tau^{*2}+ \delta^\top \widehat{\Sigma}\delta + 20\sigma_{\max}\Delta \sqrt{\frac{\log p}{n}} \norm {\delta}_1 + 2\norm {\mu_d}_2 \norm{\delta}_2.
		\end{align*}
		Plugging the above inequality into \eqref{Sigma delta infinity norm revisit}, we have
		\begin{align*}
		\norm{\hat{\Sigma}\delta}_\infty\leq \lambda\widehat{\sigma}_{\max}\left[2\tau^* +2+ \sqrt{\delta^\top \widehat{\Sigma}\delta}+\left(20\sigma_{\max}\Delta\sqrt{\frac{\log p}{n}}\norm {\delta}_1\right)^{1/2} + \left(2\norm {\mu_d}_2 \norm{\delta}_2\right)^{1/2}\right].
		\end{align*}
		Applying H{\" o}lder's inequality, we obtain that
		\begin{align*}
		\delta^\top \widehat{\Sigma}\delta \leq \lambda\widehat{\sigma}_{\max}\norm {\delta}_1\left[2\tau^* +2+ \sqrt{\delta^\top \widehat{\Sigma}\delta}+\left(20\sigma_{\max}\Delta\sqrt{\frac{\log p}{n}}\norm {\delta}_1\right)^{1/2} + \left(2\norm {\mu_d}_2 \norm{\delta}_2\right)^{1/2}\right]. 
		\end{align*}
		From the above inequality, we may derive that
		\begin{align*}
		\delta^\top \widehat{\Sigma}\delta\leq &C\cdot \lambda \sigma_{\max}\norm {\delta}_1 \Big\{\lambda\sigma_{\max}\norm {\delta}_1 + \tau^*+1+ \left(20\sigma_{\max}\Delta\sqrt{\frac{\log p}{n}}\norm {\delta}_1\right)^{1/2} \\
		&+ \left(2\norm {\mu_d}_2 \norm{\delta}_2\right)^{1/2} \Big\},
		\end{align*}
		where $C$ is a constant.
	\end{proof}

\subsection{Proof of Lemma \ref{Delta_gaps_lemma_revise}}
\begin{proof}
	When $\widehat{\tau} = \sqrt{\widehat{\beta}^\top \widehat{\Sigma}\widehat{\beta}}$, we have
	\begin{align*}
	\widehat{\tau}^2 &= \widehat{\beta}^\top \widehat{\Sigma}\widehat{\beta}
	= (\beta^* + \delta)^\top \widehat{\Sigma}(\beta^*+\delta)
	= \tau^{*2} + 2\delta^\top \widehat{\Sigma}\beta^* + \delta^\top \widehat{\Sigma}\delta\\
	&= \tau^{*2}+ \delta^\top \widehat{\Sigma}\delta +2\delta^\top (\widehat{\Sigma} - \Sigma)\beta^* + 2\delta^\top \mu_d.
	\end{align*}
	With event $\cE_1$, we have that 
	\begin{align*}
	|\widehat{\tau}^2 - \tau^{*2}|&=|\delta^\top \widehat{\Sigma}\delta +2\delta^\top (\widehat{\Sigma} - \Sigma)\beta^* + 2\delta^\top \mu_d|\\
	&\leq \delta^\top \widehat{\Sigma}\delta + 20\sigma_{\max}\Delta \sqrt{\frac{\log p}{n}} \norm {\delta}_1 + 2\norm {\mu_d}_2 \norm{\delta}_2.
	\end{align*}
	Then, using the previous results \eqref{1-norm 2-norm}, \eqref{upper_bound_revise} and  \eqref{2_norm_bound}, we obtain that
	\begin{align*}
	|\widehat{\tau}^2 - \tau^{*2}|\leq  C\cdot \Delta(\Delta+1)\sigma_{\max}^{1-q/2}M^{(3-q)/2}\sqrt{R}\left(\frac{\log p}{n} \right)^{1/2-q/4}
	\end{align*}
	for some constant $C$.
\end{proof}

\subsection{Proof of Lemma \ref{first order convergence lemma}}
\begin{proof}
Note that
\begin{align*}
\widehat{\Delta}=\sqrt{\widehat{\beta}^\top \Sigma\widehat{\beta}} &= \sqrt{\beta^\top \Sigma\beta+2\beta^\top\Sigma\delta + \delta^\top\Sigma\delta}\\
&\leq \sqrt{\beta^\top\Sigma\beta}\left(1+\frac{2\beta^\top \Sigma\delta + \delta^\top\Sigma\delta}{2\beta^\top\Sigma\beta} \right) \\
&= \Delta+\frac{2\mu_d^\top\delta+\delta^\top\Sigma\delta}{2\Delta}.
\end{align*}
Therefore, we have
\begin{align}\label{diff_1}
\frac{\Delta}{2}-\frac{\mu_d^\top\widehat{\beta}}{2\widehat{\Delta}} &= \frac{1}{2\widehat{\Delta}}(\Delta\widehat{\Delta}-\mu_d^\top\widehat{\beta})\notag\\
&\leq \frac{1}{2\widehat{\Delta}}\left(\Delta^2+\mu_d^\top (\delta-\widehat{\beta})+\frac{1}{2}\delta^\top\Sigma\delta\right)\notag\\
&=\frac{1}{4\widehat{\Delta}}\delta^\top\Sigma\delta\leq \frac{\delta^\top\Sigma\delta}{4(\Delta+\frac{\mu_d^\top \delta}{\Delta})}.
\end{align}
Note that $|\mu_d^\top\delta|\leq \norm{\mu_d}_2\norm{\delta}_2\leq M^{1/2}\Delta\norm{\delta}_2$. Using the convergence rate of $\norm{\delta}_2$ in \eqref{2-norm_bound} from Theorem \ref{1-norm_theorem}, when $n$ satisfies that
\begin{align*}
n\geq C\cdot\sigma_{\max}^2 M^{2+2/(2-q)}R^{2/(2-q)}\log p
\end{align*}
for some constant $C$, we have that $|\mu_d^\top\delta|\leq\Delta^2/2$, and thus it follows from \eqref{diff_1} that
%
\begin{align*}
\frac{\Delta}{2}-\frac{\mu_d^\top\widehat{\beta}}{2\widehat{\Delta}}\leq \frac{\delta^\top\Sigma\delta}{2\Delta}\leq \frac{M}{2\Delta}\norm{\delta}_2^2.
\end{align*}
\end{proof}

\subsection{Proof of Lemma \ref{second order convergence lemma}}
\begin{proof}
	We first show that 
	\begin{align}\label{mu hat convergence}
	\mathbb{P}\left(\norm{\widehat{\mu}^{(\ell)} - \mu^{(\ell)}}_\infty \leq \sigma_{\max}\sqrt{\frac{2\log p}{n}},~\ell=0,1  \right)\geq 1-4p^{-1}.
	\end{align}
		Note that $\widehat{\mu}^{(\ell)}\sim N(\mu^{(\ell)}, \Sigma/n)$ for $\ell = 0,1$, and thus $\widehat{\mu}^{(\ell)}_j\sim N(\mu^{(\ell)}_j, \Sigma_{j,j}/n)$, for $j\in[p]$. Hence, 
	\begin{equation*}
	\mathbb{P}\left(|\widehat{\mu}^{(\ell)}_j - \mu^{(\ell)}_j|\geq t\right)\leq 2\exp\left(-\frac{nt^2}{\Sigma_{j,j}}\right)\leq 2\exp\left(-\frac{nt^2}{\sigma^2_{\max}}\right)\quad \text{for~all~}\ell \in\{0,1\},~j\in[p].
	\end{equation*}
	Taking $t = \sigma_{\max}\sqrt{2\log p/n}$ and applying the union bound for all $j\in[p]$, we have
	\begin{equation*}
	\mathbb{P}\left(\norm{\widehat{\mu}^{(\ell)} -\mu^{(\ell)}}_\infty\leq \sigma_{\max}\sqrt{\frac{2\log p}{n}},~\ell=0,1  \right)\geq 1-4p\exp(-2\log p) = 1-4p^{-1}.
	\end{equation*}
	We next bound the term $\widehat{\beta}^\top(\mu^{(\ell)} - \widehat{\mu}^{(\ell)})$ for $\ell = 0,1$. Note that
	\begin{align*}
	\widehat{\beta}^\top(\mu^{(\ell)} - \widehat{\mu}^{(\ell)}) &= (\beta^*+\delta)^\top(\mu^{(\ell)} - \widehat{\mu}^{(\ell)}) \\
	&\leq (\norm {\beta_{S_\eta}^*}_1 +\norm {\beta_{S_\eta^c}^*}_1 + \norm {\delta}_1)\norm {\mu^{(\ell)} - \widehat{\mu}^{(\ell)}}_\infty\\
	&\leq \left(\sqrt{|S_\eta|}\norm {\beta^*}_2 +5 \norm {\beta_{S_\eta^c}^*}_1 + 4\sqrt{|S_\eta|}\norm {\delta}_2\right)\norm {\mu^{(\ell)} - \widehat{\mu}^{(\ell)}}_\infty.
	\end{align*}
	Here the last inequality uses \eqref{1-norm 2-norm}.
	Also, note that $\norm {\beta^*}\leq M^{1/2}\Delta$. With our choice of $\eta$ in \eqref{eta_choice} and the upper bound for $\norm {\delta}_2$, when $n$ satisfies that 
	\begin{align*}
	n\geq C\cdot \sigma_{\max}^2M \Delta^{-\frac{4}{2-q}}R^{\frac{2}{2-q}}\log p
	\end{align*}
	for some constant $C$, 
	we have that 
	\begin{align}\label{beta_hat mu}
	\widehat{\beta}^\top(\mu^{(\ell)} - \widehat{\mu}^{(\ell)})\leq \sigma_{\max}^{-q/2}M^{\frac{1-q}{2}} \Delta\sqrt{R}\left(\frac{\log p}{n}\right)^{\frac{2-q}{4}}.
	\end{align}
	We then consider the term $\widehat{\Delta} = \sqrt{\widehat{\beta}^\top \Sigma\widehat{\beta}}$. Note that 
	\begin{align*}
	\widehat{\Delta}^2 = \widehat{\beta}^\top \Sigma\widehat{\beta} = \Delta^2 + 2\mu_d^\top \delta + \delta^\top \Sigma\delta,
	\end{align*}
	Hence we have
	\begin{align*}
	|\widehat{\Delta}^2 - \Delta^2|\leq 2\norm{\mu_d}\norm {\delta}_2 + M\norm {\delta}_2^2.
	\end{align*}
	When $n$ is sufficiently large, we have that $|\widehat{\Delta}^2 - \Delta^2|\leq \frac{1}{2}\Delta^2$. Combining this with \eqref{beta_hat mu}, we have that  
	\begin{align*}
	\left(\frac{\widehat{\beta}^\top (\mu^{(0)} - \widehat{\mu}^{(0)}) + \widehat{\beta}^\top (\mu^{(1)} - \widehat{\mu}^{(1)})}{2\widehat{\Delta}}\right)^2
	\leq C\cdot\sigma_{\max}^{-q}M^{1-q}R\left(\frac{\log p}{n}\right)^{1-q/2}
	\end{align*}
	for some constant $C$.
	Therefore, Lemma \ref{second order convergence lemma} holds true.
\end{proof}

\section{Review of  Gautier's method}\label{review of gautier}
In this section, we provide a brief review of Gautier's pivotal method for high-dimensional linear regression in \cite{gautier2011high} that inspires our work.
Note that they  consider a more complicated high-dimensional instrumental variables model. Here we discuss the particular case where the regressors and instruments are identical for ease of presentation. Specifically, let $X\in\RR^{n\times p}$ be a design matrix with $n$ observations and $p$ variables, and let $y\in\RR^{n}$ be the response vector. We consider the following linear  model that
\begin{align*}
y = X\beta^* + \varepsilon\quad\textrm{with}\quad\varepsilon \sim N(0,\sigma^2 I_n),
\end{align*}
where $\beta^*\in \mathbb{R}^p$ is the unknown regression coefficient with $\norm{\beta^*}_0=s< n\ll p$, and $\varepsilon$ is the noise. The Gautier's estimator can be viewed as a variant of the Dantzig selector \citep{candes2007dantzig}, and is the optimal solution to the following convex optimization problem that
\begin{align}\label{self_tuned}
(\hat{\beta},\hat{\gamma}) = \argmin_{\beta, \gamma}\ &\norm{\beta}_1 + c\gamma,\\
\text{subject~to~}\ &\frac{1}{n} \norm{ X^\top(Y-X\beta)}_\infty \leq\lambda\gamma,\quad \frac{1}{n} \norm{Y-X\beta}_2^2\leq \gamma^2,\notag
\end{align}
where $c$ and $\lambda$ are two tuning parameters, and $\hat{\gamma}$ is  an estimator of $\sigma$. The theoretical analysis in \cite{gautier2011high} suggests that the tuning parameter $c$ can  be set as a constant between $0$ and~$1$, and the tuning parameter $\lambda$  can be chosen as 
\begin{align*}
\lambda = A\cdot \sqrt{\frac{2\log p}{n}},
\end{align*}
where $A$ is a constant independent of $\sigma$. 
Therefore, the Gautier's estimator is less sensitive to the parameter tuning  than the Dantzig selector, where the tuning parameter depends on $\sigma$.

\section{Numerical study on performance of Lasso, Dantzig Seector and Gautier's method}
In this section, we provide additional numerical results to compare the performance of Lasso, Dantzig Selector and Gautier's method for linear regression in high dimensions.

\begin{figure}[htb!]	
	\centering
		\includegraphics[width=0.98\textwidth]{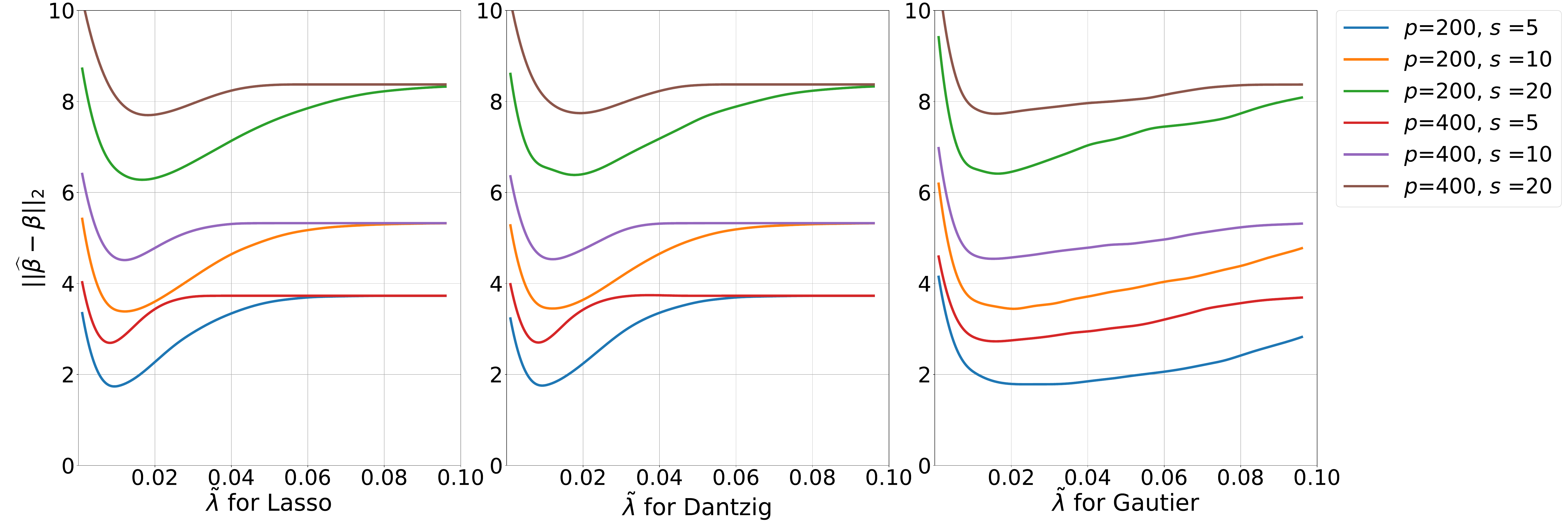}
	\caption{\it $\ell_2$ estimation error v.s. value of tuning parameter $\tilde{\lambda}$ in Lasso (left), Dantzig selector (middle) and Gautier's method (right).
		The results are averaged over 100 replicates.}
	\label{regression_error_2}
\end{figure}

We generate the data by a process considered in \cite{candes2007dantzig}. To be more specific, we set $n=100$, $p=200$, $s=5,10,20$. We generate the rows of $X$ from the standard Gaussian distribution and then normalize each row of $X$. 
For $\beta^*$, we set
\begin{align*}
	\beta^*_i = u_i(1+|a_i|) \text{~for~}i=1,\cdots,s,
\end{align*}
where $u_i=\pm1$ with probability 1/2, and $a_i\sim N(0,1)$ and is independent of $u_i$. Meanwhile, we set $\sigma = \sqrt{\frac{s}{n}}$.
To fine-tune the parameter, we generate an independent validation set with same sample size $n=100$ as the training set. We let $\lambda = \tilde{\lambda}\sqrt{\frac{\log p}{n}}$ for all the three methods, and we tune the factor $\tilde{\lambda}$ over a range from 0 to 1 for each method. Figure \ref{regression_error_2} shows the results of the estimation error  $\|\hat\beta-\beta^*\|_2$
versus the $\tilde{\lambda}$ value in the three methods, averaged over 100 replicates under each setting of different $p$ and $s$. 
For Gautier's method, the result is not sensitive to the parameter $c$ as long as $c$ is not too small, and we set $c=20$. Table \ref{table: regression} summarizes the  estimation error $\|\hat\beta-\beta^*\|_2$ under different $p$ and $s$.
As can be seen, the three methods have comparable performance in $\beta^*$ estimation after fine-tuning.
\begin{table}[!htb]	
	\centering
	\caption{\it The $\ell_2$ error of $\beta^*$ estimation for the regression example. The testing errors are averaged over 100 replicates. The standard deviation of the testing errors are given in brackets. }
\vspace{0.1in}
	\begin{tabular}{p{3cm}<{\centering}|cccccc}
		\hline
		\multirow{2}*{Method}& \multicolumn{6}{c}{$(s,p)$} \\
		&$(5, 200)$&$(10, 200)$&$(20,200)$&$(5, 400)$&$(10, 400)$&$(20, 400)$\\
		\hline
		\multirow{2}*{Lasso}&1.801& 3.425 &6.546&2.738&4.609&7.748\\
		&(0.325)&(0.517)&(0.969)&(0.440)&(0.492)&(0.460)\\
		\hline
		\multirow{2}*{Dantzig Selector}&1.802& 3.466&6.389&2.757&4.600&7.744\\
		&(0.343)&(0.495)&(0.634)&(0.449)&(0.500)&(0.471)	\\
		\hline
		\multirow{2}*{Gautier's Method}&1.771& 3.412&6.375&2.749&4.653 & 7.741\\
		&(0.341)&(0.406)&(0.620)& (0.441) &(0.548)&(0.483)	\\
		\hline
	\end{tabular}
\label{table: regression}
\end{table}

\section{Technical derivation on the penalty term in PANDA}\label{penalty_term}

In this section,  we provide a deep insight on how to non-trivially modify Gautier's pivotal method to our context. To be more specific, we compare the penalty term imposed in Gautier's pivotal method and our proposed PANDA, and explain our choice of a quadratic penalty for $\tau$ in (\ref{con:problem}). For simplicity, we consider the case where $q=0$ and $|\supp(\beta^*)|\leq s$.

Let $S=\supp(\beta^*)$. For both Gautier's method and PANDA, a key step to derive the upper bound of $\norm{\delta}_1 = \norm{\widehat{\beta}-\beta^*}_1$ is to show that $\delta$ belongs to some restricted subset $\cC_{S,\beta^*}$ with high probability, where $\cC_{S,\beta^*}$ is defined in (\ref{cone definition}). Note that when $q=0$, $\norm{\beta^*_{S^c}}_1 =0$, such that $\cC_{S,\beta^*}$ reduces to
\begin{align*}
\cC_S = \left\{\delta\in \mathbb{R}^p:\norm{\delta_{S^c}}_1\leq 3\norm{\delta_S}_1 \right\}.
\end{align*}

In Gautier's method, it is shown that with high probability, $(\beta^*,\sigma^*)$ is feasible to the program (\ref{self_tuned}), where $\beta^*$ is the true regression parameter and
\begin{align*}
\sigma^* \coloneqq \frac{1}{\sqrt n}\norm{Y-X\beta^*}_2.
\end{align*}
Then, by the optimality condition of the solution $\widehat{\beta}$, i.e.
$\norm{\widehat{\beta}}_1 + c\widehat{\sigma}\leq \norm{\beta^*}_1 + c\sigma^*$, $\norm{\delta_{S^c}}_1$ can be upper bounded by
\begin{align*}
\norm{\delta_{S^c}}_1&\leq \norm{\delta_S}_1 + \frac{c}{\sqrt{n}}\left(\norm{Y-X\beta^*}_2 -
\norm{Y-X\widehat{\beta}}_2 \right)\\
&\leq \norm{\delta_S}_1 + \frac{c}{\sqrt{n}}\delta^\top
\frac{X^\top (Y-X\beta^*)}{\norm{Y-X\beta^*}_2}\\
&\leq \norm{\delta_S}_1 + c\norm{\delta}_1\frac{\norm{\frac{1}{n}X^\top(Y-X\beta^*)}_\infty}{\sigma^*}\\
&\leq \norm{\delta_S}_1 + c\lambda\norm{\delta}_1,
\end{align*}
where the second inequality uses the convexity of $\norm {Y-X\beta}_2$ in $\beta$, the third inequality uses H{\" o}lder's inequality and the definition of $\sigma^*$, and the last inequality is due to the first constraint in (\ref{self_tuned}). With properly chosen $c$ and $\lambda$, 
it can be shown that 
$\delta\in\cC_S$ with high probability.

For PANDA, if we follow the above framework and impose the same penalty $c\tau$, a similar argument leads to
\begin{align*}
\norm {\delta_{S^c}}_1&\leq \norm {\delta_S}_1 + c\norm{\delta}_1\frac{\norm {\widehat {\Sigma}\beta^*}_\infty}{\sqrt{\beta^*\widehat{\Sigma}\beta^*}}\\
&\leq \norm {\delta_S}_1 + c\norm{\delta}_1\frac{\norm {\widehat {\Sigma}\beta^* - \widehat{\mu}_d}_\infty + \norm{\widehat{\mu}_d}_\infty}{\sqrt{\beta^*\widehat{\Sigma}\beta^*}}\\
&\leq \norm {\delta_S}_1 + c\norm{\delta}_1 \left(\lambda+ \frac{\lambda + \norm{\widehat{\mu}_d}_\infty}{\sqrt{\beta^{*\top}\widehat{\Sigma}\beta^*}}\right).
\end{align*}
Note that $\sqrt{\beta^{*\top}\widehat{\Sigma}\beta^*}$ converges to $\Delta$, and thus the term $\frac{\norm{\hat{\mu}_d}_\infty}{\Delta}$ dominates the last term, and the choice of $c$ must rely on the unknown $\Delta$ to ensure that $\delta\in\cC_S$ with high probability. 

In other words, we cannot directly follow Gautier's framework to impose the penalty $c\tau$. Nevertheless, Gautier's method inspires us to impose a quadratic penalty term on $\tau$, by which it turns out that the tuning parameters will no longer rely on the unknown $\Delta$.

Here we remark that in order to guarantee the tuning-insensitive property of our PANDA method, the penalty on $\tau$ must be quadratic.
Suppose we consider an increasing and convex penalty function $f(\tau)$ instead. Technically, in order to guarantee that $\delta = \hat{\beta} - \beta^*$ belongs to the restricted set 
\begin{align*}
	\cC_{S,\beta^*} \coloneqq \left\{\delta\in\mathbb{R}^p:~ \norm{\delta_{S^c}}_1\leq 3\norm{\delta_S}_1 + 4\norm{\beta^*_{S^c}}_1 \right\}
\end{align*}
with high probability, we require 
$f$ to satisfy that $f(\tau^*) - f(\hat{\tau})\leq  \frac{1}{2}\norm{\delta}_1$, where $\tau^*$ is close to $\Delta$.
Following the argument in the proof of Lemma 7, we can derive an upper bound for $f(\tau^*) - f(\hat{\tau})$ as follows:
\begin{align*}
	f(\tau^*) - f(\hat{\tau})&\leq f\left(\sqrt{\beta^{*\top}\hat{\Sigma}\beta^*}\right) - f\left(\sqrt{\hat{\beta}^\top\hat{\Sigma}\hat{\beta}}\right)\\
	&\leq\left|\frac{f^\prime(\tau^*)}{\tau^*}\right|\norm{\hat{\Sigma}\beta^*}_\infty\norm{\delta}_1\\
	&\leq \left|\frac{f^\prime(\tau^*)}{\tau^*}\right|\left(\norm{\mu_d}_\infty+\norm{(\hat{\Sigma} - \Sigma)\beta^*}_\infty\right)\norm{\delta}_1.
\end{align*}
In order to control that $f(\tau^*) - f(\hat{\tau})\leq \frac{1}{2}\norm{\delta}_1$, we need $\left|\frac{f^\prime(\tau^*)}{\tau^*}\right|\left(\norm{\mu_d}_\infty+\norm{(\hat{\Sigma} - \Sigma)\beta^*}_\infty\right)\leq \frac{1}{2}$. When $n$ is sufficiently large, the term $\norm{(\hat{\Sigma} - \Sigma)\beta^*}_\infty$ here is small, and $\norm{\mu_d}_\infty$ can be closely estimated  from the sample. 
Therefore, we require the term $\frac{f^\prime(\tau^*)}{\tau^*}$ to be controlled by some constant that is independent of $\tau^*$ or $\Delta$. To satisfy this, the Taylor expansion of $f$ can only have non-zero coefficient for the first-order term, while the coefficients for other orders must be zero, implying that $f$ is a quadratic function.

\vskip 0.2in
\bibliography{sample}

\end{document}